\documentclass[matzei,final,envcountsame]{svjour}
\usepackage{amssymb,euscript}
\usepackage[all]{xy}

\newcommand{\adel}{{\cA}_\partial}
\newcommand{\cA}{{\mathcal A}}
\newcommand{\cE}{{\mathcal E}}
\newcommand{\cF}{{\mathcal F}}
\newcommand{\cG}{{\mathcal G}}
\newcommand{\cin}{C^\infty}
\newcommand{\cN}{{\mathcal N}}
\newcommand{\cS}{{\EuScript S}}
\newcommand{\cV}{{\mathcal V}}
\newcommand{\dem}{{\noindent{\bf {Proof: }}}}
\newcommand{\acV}{{^a\! \cV}}
\newcommand{\assN}{{^a\! S^*\! N}}
\newcommand{\atm}{{^a\! TM}}
\newcommand{\atmo}{\atm_{|t=0}}
\newcommand{\atx}{{^a\! TX}}
\newcommand{\atsM}{{{}^a\!T^* M}}
\newcommand{\atsx}{{{}^a\!T^* X}}

\newcommand{\cz}{{\mathbb C}}
\newcommand{\End}{\mathrm{End}}
\newcommand{\ff}{\mathrm{ff}}
\newcommand{\ffa}{\ff_a}
\newcommand{\hol}{{\mathrm{hol}}}
\newcommand{\nabf}{{\nabla^{b\!f}}}

\newcommand{\ind}{{\mathrm{index}}}
\newcommand{\nz}{{\mathbb N}}
\newcommand{\oeta}{\overline{\eta}}
\newcommand{\pa}{\Psi_a}

\newcommand{\ps}{{\partial_\theta}}

\newcommand{\rz}{{\mathbb R}}

\newcommand{\sfl}{{\mathrm{sf}}}
\newcommand{\Spec}{{\mathrm{Spec}}}
\newcommand{\sus}{{\mathrm{sus}}}
\newcommand{\supp}{{\mathrm{supp}}}
\newcommand{\tna}{{\tilde{\nabla}}}
\newcommand{\topp}{{\mathrm{top}}}
\newcommand{\tps}{{\tilde\partial_\theta}}
\newcommand{\tr}{{\mathrm{tr}}}
\newcommand{\Tr}{{\mathrm{Tr}}}
\newcommand{\zz}{{\mathbb Z}}

\begin{document}

\title{Adiabatic limits of eta and zeta functions 
of elliptic operators}
\author{Sergiu Moroianu\thanks{Partially supported by
ANSTI (Romania), the European
Commission RTN HPRN-CT-1999-00118 \emph{Geometric Analysis} and by the IREX RTR
project.}}
\date{\today}
\institute{Institutul de Matematic\u{a} al Academiei Rom\^{a}ne, 
P.O. Box 1-764,\\ RO-70700
Bucharest, Romania\\
\email{moroianu@alum.mit.edu}}
\maketitle
\begin{abstract}
We extend the calculus of adiabatic pseudo-differential operators
to study the adiabatic limit behavior of the eta and zeta functions
of a differential operator $\delta$, 
constructed from an elliptic family of operators indexed by $S^1$. 
We show that the regularized values ${\eta}(\delta_t,0)$ and 
$t{\zeta}(\delta_t,0)$ are smooth functions of $t$ at $t=0$, and we 
identify their values at $t=0$ with the holonomy of the determinant bundle, 
respectively with a residue trace.
For invertible families of operators,
the functions ${\eta}(\delta_t,s)$ and $t{\zeta}(\delta_t,s)$ 
are shown to extend smoothly to $t=0$ for all values of $s$. 
After normalizing with a Gamma
factor, the zeta function satisfies in the adiabatic limit an identity
reminiscent of the Riemann zeta function, while the eta function converges to
the volume of the Bismut-Freed meromorphic family of connection $1$-forms.

\end{abstract}

\vspace{0.5cm}\section{Introduction}

The eigenvalues of elliptic operators on compact manifolds possess miraculous
properties.
Starting with the foundational work of Minakshisundaram and Pleijel
\cite{mipl}, we know for instance that the zeta function 
$$s\mapsto \sum_{\Spec(D)\ni\lambda\neq 0} |\lambda|^{-s/2}$$
of a second order self-adjoint elliptic pseudo-differential operator $D$
on a closed manifold $M$ 
is analytic in the half-plane $\{\Re(s)>\dim(M)\}$ and extends analytically to $\cz$ with simple 
poles at $\dim(M)-\nz$.
If $D$ is self-adjoint of order 1, a refinement of
the zeta function of $D^2$ is the eta function of Atiyah, Patodi and Singer \cite{aps},
$$s\mapsto \sum_{\Spec(D)\ni\lambda\neq 0}
\mathrm{sign}(\lambda)|\lambda|^{-s}.$$ 
This function is also meromorphic on $\cz$ with simple poles at $\dim(M)-\nz$. 

To grasp the significance of these facts, one must take into account that it is
generally impossible to determine all the eigenvalues of a given operator, or to
decide which discrete subsets of $\rz$ can occur as spectra of self-adjoint 
elliptic operators. The pole at $s=\dim(M)$ of the zeta function is related to the Weyl
asymptotic formula for the eigenvalues, but the next poles are more subtle.
Even when one can find explicitly the eigenvalues, there is no obvious reason why 
the series defining the zeta function should extend analytically to $\cz$. 
In the simplest non-trivial case where $M$ is the circle and
$D:=i\ps$ is the angular derivative, the zeta function of
the Laplacian $D^2$ equals the
Riemann zeta function up to a factor of $2$. Beyond the usual elementary tricks 
there exists
therefore a more fundamental reason for the analytic extension of this classical
function, that is the fact that $\zz$ is the spectrum of a differential
operator. 

One approach to the analysis of $\zeta(s)$
particularly unsuitable for generalizations seems to be the 
infinite product formula over the primes. Number theorists introduce 
in this product an additional factor $\Gamma(s/2)$, corresponding to the "prime" $0$. 
Then the normalized Riemann zeta function $\overline{\zeta}(s):=\Gamma(s/2)
\zeta(s)$ satisfies the functional equation
$$\overline{\zeta}(1-s)=\pi^{1/2-s}\overline{\zeta}(s).$$

Amazingly enough, we have found in the context of 
adiabatic limits a functional identity (Theorem \ref{mz}) 
in terms of similarly normalized zeta 
functions. In this context, the adiabatic limit of the eta function
appears to be related to the meromorphic family of connection $1$-forms
on the determinant line bundle constructed by Bismut and Freed \cite{bf2} if,
again, both are
renormalized with appropriate Gamma factors (Theorem \ref{eil}). 

Our initial motivation was to give a simple proof of the holonomy formula of
\cite{bf2} relating the adiabatic limit of the eta invariant to the
determinant bundle of a family of Dirac operators on the circle as conjectured by
Witten \cite{wit}. However, the results turn out to be valid for the normalized 
eta and zeta
functions themselves and not just for their values at $s=0$ (for comparison, while it
is interesting to know that the Riemann zeta function equals $-1/2$ at $s=0$, it is
certainly desirable to say something about it at other points as well).
As a corollary we give a formula for the adiabatic limit of the
determinant.

This paper treats general elliptic first order differential operators and not only
Dirac operators, which is in our view a significant breakthrough. 
Indeed, previous works on eta invariants use the geometric construction of the heat 
kernel and the so-called rescaling technique of Getzler, which is 
intimately related to the Clifford algebra. These methods are unlikely to extend 
to the general case, therefore suggesting that Dirac operators enjoy 
special properties among first-order
operators. We show that this is not the case for adiabatic
limits of the eta and zeta functions (albeit it is certainly true for statements like
the local index theorem).

Our approach is therefore completely different. We construct the calculus
$\Psi_{ae}(X;M)$ of extended adiabatic operators, which extends 
(as the notation shows) the calculus $\pa(N)$ of adiabatic 
pseudo-differential operators introduced in \cite{mame90}, 
see also \cite{nwx} for an approach using differentiable groupoids. 
The calculus $\pa(N)$ was originally constructed for the study
of adiabatic limits in the framework of Melrose's programme for quantizing 
singular geometric structures \cite{mel1}. An adiabatic limit means blowing up the 
metric on the base of a fibration of manifolds by a factor $t^{-2}$ as $t\to 0$. 
The Laplacian of this metric will stop being elliptic when the 
parameter $t$ reaches $0$. The basic
idea is to force this operator to be elliptic at $t=0$ in an appropriate sense.
Instead of heat operators we use complex powers, which are straightforward to
construct inside the adiabatic algebra.

A immediate gain of our point of view is replacing the limit as
$t\to 0$ by a 
Laurent-type singularity at the end point of the closed interval
$[0,\infty)$. In particular, we 
deduce, under some invertibility hypothesis, that the eta function 
has a Taylor expansion in the adiabatic limit, thus we can rule out 
for instance $\log t$-like terms on a priori grounds. 

The regularity at $s=0$ of the eta function is not used in the proof. This is
actually not so surprising, since our formulas are valid (in the sense of Laurent 
coefficients) also at the poles of the eta function.

Let us describe briefly the contents of the paper. In Section \ref{prelim}
we introduce the basic objects and give an overview of the results. 
The Dirac operator on the total space of a fibration over the circle 
is analyzed in Section \ref{dir}.
The determinant line bundle of families over $S^1$ 
is reviewed in Section \ref{det}, and analytic extensions
of zeta-type functions in Section \ref{aert}. Section \ref{ea} deals with the 
adiabatic algebra and its properties. The main results are stated and proved in 
Section \ref{mate} in the invertible case and in Section \ref{gen} in the general case. 
Finally, in Section \ref{clorem} we state without proof some related results.

The results were announced in a note in Comptes rendus \cite{cr}.

\begin{acknowledgement} I owe to Richard Melrose the idea of applying the 
adiabatic algebra to the study of the eta invariant. Jean-Michel Bismut 
showed to me how to treat the non-invertible case in Section \ref{gen}. 
Section \ref{dir} is due to Andrei Moroianu (although Proposition \ref{am} 
may be extracted from \cite{bf2,bc}). The appearance of Gamma
functions in the formulas was suggested by Nicu\c{s}or Dan. I have 
also benefited from the advice of Christian B\"{a}r, Mattias Dahl and
Robert Lauter. Last but not least, I am grateful to the anonymous
referee for very pertinent remarks and for
detecting many mistakes, both in presentation 
and in substance.
\end{acknowledgement}

\vspace{0.5cm}\section{Preliminaries}\label{prelim}

Let $N$ be a closed Riemannian manifold, $E$ a Hermitian vector bundle
over $N$ and $\delta$ a first order elliptic differential operator on 
the sections of $E$. We define the normalized zeta function of $\delta$ by 
\begin{eqnarray}\label{zeta}
\overline{\zeta}(\delta,s)&:=&\Gamma\left(\frac{s}{2}\right)
\sum_{0\neq\lambda\in\Spec(\delta^*\delta)}
\lambda^{-s/2}+\dim\ker(\delta).
\end{eqnarray}
If $\delta$ is self-adjoint, we define similarly its normalized eta function by
\begin{eqnarray*}\label{eta}
\oeta(\delta,s)&:=&
\pi^{-\frac{1}{2}}\Gamma\left(\frac{1+s}{2}\right)
\sum_{0\neq\lambda\in\Spec(\delta)}\mathrm{sign}
(\lambda)|\lambda|^{-s}
+\dim\ker(\delta).
\end{eqnarray*}
These functions are holomorphic in the half-plane $\{\Re(s)>\dim N\}$
and extend analytically to $\cz$ with possible simple poles at
$\dim N-2k$, respectively $\dim N-2k-1$, for $k\in\nz$ (Corollary \ref{cme}). 
The definitions differ from the usual zeta and eta functions by some Gamma
factors which will turn out to be the key for the validity of our 
results outside $s=0$.
By the result of \cite{aps}, ${\eta}(\delta,s)$ is finite at $s=0$,
hence the regularized value of $\oeta(\delta,s)$ at $s=0$ coincides with the
refined eta-invariant of Atiyah-Patodi-Singer. 
We stress however that the results of this paper are independent of this
regularity. 

Suppose that $N$ is the total space of a locally trivial fibration
of compact manifolds $N\stackrel{p}{\longrightarrow}M$. We are
interested in the restriction of this fibration to a circle in the 
base, we can therefore assume that $M=S^1$. Let $E=E^+\oplus E^-$ be a 
$\zz/2\zz$-graded Hermitian bundle over $N$ and $D$ a family of 
elliptic differential operators of order 1 over the fibers of $p$.
In other words, for every $x\in S^1$, the operator $D_x:\cin(N_x,E^+|_{N_x})
\rightarrow \cin(N_x,E^-|_{N_x})$ is elliptic and depends smoothly on 
$x\in S^1$. 

We fix a connection in the fibration $N\rightarrow S^1$ and a family 
of metrics $g$ on the fibers extended by $0$ on the horizontal 
distribution. For any vector $Y$ tangent to $S^1$ we denote by 
$\tilde Y$ its horizontal lift to $N$. We also fix differential 
operators $\nabla_{\tilde{Y}}$ on the sections of $E^\pm$ with the 
property that 
$(\nabla_{\tilde{Y}}e_1,e_2)+(e_1,\nabla_{\tilde{Y}}e_2)=\tilde{Y}(e_1,e_2)$. 
Such operators arise for instance as the restriction to horizontal 
vectors of metric connections in $E^\pm$. Let $\tr(L_{\tilde{Y}}g)$ be 
the contraction by $g$ of the Lie derivative of the tensor $g$. Since 
$g$ is vertical and $\tilde{Y}$ is horizontal, this expression is
tensorial in $Y$. One sees easily that the operators
$$\tna_Y:=\nabla_{\tilde{Y}}+\frac{1}{4}\tr(L_{\tilde{Y}}g)$$
are skew-symmetric. If $|d g|$ is the volume density induced by the metric $g$
then $L_{\tilde{Y}}|d g|=\frac{1}{2}\tr(L_{\tilde{Y}}g)|d g|$.

On $S^1$ we fix the canonical metric $d\theta^2$ (any other metric is 
isometric to a multiple of $d\theta^2$). Let $\ps$ be the 
positively oriented unit vector field on $S^1$. For $t\in[0,\infty)$
consider the self-adjoint operator
\begin{equation}\label{deltat}
\delta_t:=\left[
\begin{array}{cc}
-ti\tna_{\ps}& D^*\\ D& 
ti\tna_{\ps}
\end{array}
\right]. 
\end{equation}
acting on $\cin(N,E)$.
Such an operator arises for instance in \cite{bf2},
where the fibers of $N\rightarrow S^1$ are even-dimensional and carry a continuous
spin structure, $E^\pm$ are the spinor bundles of the fibers and $\nabla$ 
is the Levi-Civita connection. In that case $\delta_t$ is the Dirac operator on $N$
associated to the metric
\begin{equation}\label{gt}
g^N_t:=g+\frac{d\theta^2}{t^2}.
\end{equation}
Note that $\delta_t$ is elliptic if and only if $t>0$;
nonetheless, $\delta_t$ is an elliptic \emph{adiabatic} family of operators 
including at the
limiting value $t=0$, in the sense of Section \ref{ea}. 

The purpose of this paper is to relate spectral invariants of the family 
$\{D_x\}_{x\in S^1}$ on the fibers to those of the operators 
$\{\delta_t\}_{t\in(0,\infty)}$. 

In the case when $D$ is a family of twisted Dirac operators, a conjecture of Witten 
\cite{wit} proved by Bismut and Freed \cite{bf2} states that
$\lim_{t\rightarrow 0} \exp({-i\pi\eta(\delta_t)})$
exists and equals the holonomy of the determinant line bundle $\det(D)$ over $S^1$.
This turns out to be true for any family $D$ of first order elliptic differential
operators, and moreover we prove that
$\oeta(\delta_t)$ is smooth in $t$ (modulo $2\zz$) for $t\in[0,\infty)$. 
If the family $D$ is
invertible, we prove that the function $\oeta(\delta_t,s)$ admits a Taylor
expansion at $t=0$ with coefficients meromorphic functions on $\cz$.
The limit as $t\to 0$ is computed explictly in terms of the connection $1$-form on the
determinant line bundle.

Similar results hold for the zeta function. In the invertible case, we show that 
the family $t\overline{\zeta}(\delta_t,s)$ admits a Taylor expansion in powers of 
$t$ with coefficients meromorphic functions on $\cz$. We show that
$\lim_{t\rightarrow
0}t\overline{\zeta}(\delta_t,s)$ equals the average of $\overline{\zeta}(D,s-1)$
over $S^1$.
As corollaries we obtain formulas for the adiabatic 
limit of $t\zeta(\delta_t,0)$ and of $t\log\det(\delta_t)$. The later formula 
takes a simpler form when the manifold $N$ is odd-dimensional. The former is 
valid for general elliptic families $D$.

We close this section with a remark about the Gamma functions. The adiabatic 
limit process involves passing from
a $(n+1)$-dimensional manifold to a family of $n$-dimensional manifolds. 
From zeta-type functions of operators on $N$ one obtains in the limit a factor 
$$
f(s):=\int_{-\infty}^\infty (1+\tau^2)^{-s}d\tau
$$
where $\tau$ encodes the missing dimension. We note the identity
\begin{equation}\label{fs}
f(s)=\frac{\sqrt{\pi}\Gamma\left(s-\frac{1}{2}\right)}{\Gamma(s)}.
\end{equation}

\vspace{0.5cm}\section{The lifted Dirac operator}\label{dir}

This section computes an important example of lifting an operator from the fibers 
of a fibration to the total space. This example motivates the definition
(\ref{deltat}) of the family $\delta_t$. 

Let $p:X\rightarrow M$ be a fibration of closed
manifolds, whose fibers are oriented Spin Riemannian manifolds of
dimension $2n$, and assume that there exists a continuous spin 
structure on the fibers, meaning that $TX/M$ is associated to a principal 
$Spin_{2n}$ bundle which is a double cover of the orthonormal frame 
bundle of the fibers. We denote by $\{D_x\}_{x\in M}$ the
associated family of Dirac operators. We fix a connection in the fibration 
$p$, i.e.\ a choice of a horizontal distribution in $X$. This
allows us to extend the family of fiberwise metrics to a 2-tensor $g$ on
the total space of $X$. Any metric $g^M$ on $M$ induces by
pull-back, together with $g$, a Riemannian metric on $X$:
$$g^X=g+p^* g^M.$$

We define a connection $\nabla$ on the vertical sub-bundle 
$TX/M$ of $TX$ as the vertical projection of 
the Levi-Civita connection of the metric $g^X$ on $X$. 
This depends on the tensor $g$, but not on the metric $g^M$ on the base.
We also denote by $\nabla$ the induced connection on the spinors over the 
fibers. For all vectors fields $Y$ on $M$, we define a first-order 
differential operator $\tna_Y$ on sections of the 
spinor bundles by 
\begin{equation}\label{cbf}
\tna_Y:=\nabla_{\tilde{Y}}+\frac{1}{4}\tr(L_{\tilde{Y}}g).
\end{equation}
Here $\tr$ means the contraction by $g$ and $L_{\tilde{Y}}g$ is the Lie 
derivative of the tensor $g$ in the direction of the horizontal lift
$\tilde{Y}$ of $Y$. If $\alpha$ is a vertical $1$-form, then for all $f\in\cin(M)$ we
have $L_{\tilde{fY}}(\alpha)=p^*f L_{\tilde{Y}}(\alpha)$, hence $\tr(L_{\tilde{Y}}g)$
is tensorial in $Y$. Therefore, 
$\tna$ can be thought of as a connection in the
(infinite dimensional) vector bundle over $M$ of fiber-wise sections
of the spinor bundles, and the correction term ensures that this connection 
preserves the $L^2$ metric. It is $\tna$ that occurs both in the 
Bismut-Freed connection form and in the lifted Dirac operator.

Let $S^1$ be a smooth loop embedded in $M$, and
$N:=p^{-1}(S^1)$. Fix the trivial (non-bounding) spin structure on $S^1$.
Then $N$ inherits a spin structure, 
and moreover the spinor bundle is isomorphic to 
$S^+ \oplus S^-$. The Clifford action of the unit horizontal vector
$\ps$ is simply $\mp i$ on $S^\pm$. This is
one of the two possible choices, the other one being $\pm i$.
We choose this sign in order to get Witten's original sign in the
holonomy formula \cite{wit}.

Let now $\delta_t$ be the Dirac operator for the base metric scaled
by $t^{-2}$, i.e.\ for the metric
\begin{equation}
g_t^N=g+t^{-2}g^M.
\end{equation}
\begin{proposition}\label{am}
In matrix form, using the splitting $S=S^+\oplus S^-$, we have
$$\delta_t=\left[
\begin{array}{cc}
-ti\tna_{\ps}& D^*\\ D& ti\tna_{\ps}
\end{array}
\right].$$
\end{proposition}
\dem
The Levi-Civita connection $\nabla^t$ on $(N,g_t^N)$ satisfies the
identities
\begin{eqnarray*}
\nabla^t_U V &=& \nabla_U V -\frac{t^2}{2}(L_{\tps}g)(U,V) \tps\\
(\nabla^t_{t\tps}U,V)&=&t(\nabla_{\tps}U,V),
\end{eqnarray*}
where $U,V$ are vertical vectors and $\tps$ is the horizontal lift 
of the unit
tangent vector to $S^1$ for the initial metric $g^M$.
Let $\phi$ be a section of the spinor bundle $S=S^+\oplus S^-$.
Locally, $\phi$ is determined by a pair $[P,\sigma]$, where $P$ is a
local section of the Spin bundle which sits over a local
orthonormal frame $(e_1,\ldots,e_{2n},t\tps)$, and $\sigma$ lives
in the standard spinor representation $\Sigma_{2n+1}$.
For any vector $V\in TN$, we have
\begin{eqnarray*}
\nabla^t_V(\phi)&=&[P,V(\sigma)]+\frac{1}{2}\sum_{i<j\leq 2n} 
(\nabla^t_V e_i,e_j)c(e_i)c(e_j)\phi\\
&& +\frac{1}{2}\sum_{i\leq 2n} (\nabla^t_V e_i, t\tps) c(e_i) c(t\tps)\phi.
\end{eqnarray*}
Recall that $t\tps$ has length $1$ for $g_t^N$ and that $c(t\tps)$ 
equals $\mp i$ on $S^\pm$. 
Therefore
\begin{eqnarray*}
\delta_t\phi&=&\sum_{k=1}^{2n}c(e_k)\left( [P,e_k(\sigma)]
+\frac{1}{2}\sum_{i<j\leq 2n} (\nabla^t_{e_k} e_i,e_j)
c(e_i)c(e_j)\phi\right.\\ && \left. +\frac{1}{2}\sum_{i\leq 2n} 
(\nabla^t_{e_k} e_i, t\tps) c(e_i) c(t\tps)\phi\right)
+c(t\tps) \nabla^t_{t\tps}(\phi)\\
&=& \left[\begin{array}{cc}-it\nabla_{\tps}&D^*\\ D &it\nabla_{\tps}
\end{array}\right]
\phi+\frac{1}{2}\sum_{k=1}^{2n}\sum_{i\leq 2n} 
(\nabla^t_{e_k} e_i, t\tps) c(e_k) c(e_i) c(t\tps)\phi\\&=&
\left(\left[\begin{array}{cc}-it\nabla_{\tps}&D^*\\ D &it\nabla_{\tps}
\end{array}\right]
-\frac{t}{4}c(t\tps)\sum_{k,i=1}^{2n} (L_{\tps}
g)(e_k,e_i) c(e_k) c(e_i)\right)\phi\\&=&
\left(\left[\begin{array}{cc}-it\nabla_{\tps}&D^*\\ D &it\nabla_{\tps}
\end{array}\right]+\frac{t}{4}c(t\tps)\tr(L_{\tps}g)\right)\phi.
\end{eqnarray*}
\qed

\vspace{0.5cm}\section{Analytic extensions and regularized traces}\label{aert}

This section contains a review of the basic results allowing one 
to define analytic extensions of zeta-type functions. Proofs are included
for the benefit of the reader.

Recall that classical pseudo-differential operators on a closed manifold 
$X$ have a characterization
in terms of their Schwartz kernels on $X^2$. Namely, such a kernel is a 
distribution $\psi$ on $X^2$ with $sing\, supp (\psi)$ contained in the diagonal
$\Delta_X$; moreover, $\psi$ must be \emph{classical conormal} to $\Delta_X$, 
in the sense that the Fourier transform in the normal 
directions of a cut-off of $\psi$ near $\Delta_X$ is a classical symbol, i.e.\ a 
symbol admitting an asymptotic expansion in homogeneous components with step 1. 
Since we work with complex powers of operators, we allow classical symbols of any 
complex orders. By Peetre's Theorem, the differential operators are exactly 
those pseudo-differential operators with Schwartz kernel supported on the diagonal. 

\begin{proposition}\label{tilp}
Let $P(s)$ be an entire family of classical pseudo-differential operators on
$X$, such that $P(s)\in\Psi^s(X)$. Then ${\Tr}(P(s))$, which is 
\emph{a priori} defined for $\Re(s)<-\dim X$, extends analytically to
$\cz$ with possible simple poles at integers $k\geq-\dim X$. \label{it2}
\end{proposition}
\dem For $\Re(s)<-\dim X$, $\Tr(P(s)$ equals the integral of the pointwise 
trace of the distributional kernel of $P(s)$ over the diagonal. By pulling back to the 
tangent bundle and then Fourier transforming in the fibers, this becomes
the integral on $T^* X$ of the pointwise trace of the
full symbol of $P(s)$, which is an entire family $\{a(s)\}_{s\in \cz}$ such that 
$a(s)$ is a classical symbol of order $s$. The claim follows from the next Lemma.
\qed

\begin{lemma}\label{bl}
Let $a(s)$ be an entire family of symbols of order $s$ on a vector bundle $V$
over $X$. Fix a density $d\mu$ on $X$ and a Hermitian 
metric $g$ in $V$. Then 
$\int_V a(s)d\mu dg$, which is well-defined for $\Re(s)<-\dim V$,  
extends analytically to $s\in \cz$, with at most 
simple poles at real integers.
\end{lemma}
\dem If $a_k(s,v)$ is an entire family of symbols on $V$ such that
for $\|v\|\geq 1$, $a_k(s,v)=h\left(s,\frac{v}{\|v\|}\right)\|v\|^{k+s}$ where 
$h$ is entire in $s$, then 
$\int_V a(s,v)dv$ splits in the integral on the unit ball bundle (which is entire) 
and the integral on the exterior of the unit ball, which can be computed 
explicitly in polar coordinates. We denote by $S(V)$ the sphere bundle inside
$V$.
\begin{eqnarray*}
\int_{|v|\geq 1}
a(s,v)dv&=&\int_1^\infty\int_{S(V)}h(s,\theta)r^{k+s+\dim V-1}d\theta dr\\
&=&-\frac1{s+k+\dim V}\int_{S(V)}h(s,\theta)d\theta.
\end{eqnarray*} 
Any family of classical symbols $a(s)$ 
of order $s$ can be decomposed as $a(s)=a_0(s)+\tilde{a}(s)$, where $a_0(s,v)$
is a family as above for $k=0$, and $\tilde{a}(s)$ is an entire family of 
symbols of order $s-1$. Then $\int a_0(s,v)dv$ has just one pole at 
$s=-\dim V$, while $\int \tilde{a}(s,v)dv$ is well-defined for 
$\Re(s)<-\dim V+1$. The Lemma follows by iteration of this argument.
\qed

Recall the definition of the Wodzicki residue trace: fix an elliptic 
first-order positive
pseudo-differential operator $Q$ on the $m$-dimensional manifold $X$. 
Then the complex powers of $Q$ form an entire family of classical
pseudo-differential operators \cite{seeley}.
By Proposition \ref{tilp}, for any classical pseudo-differential operator 
$P\in\Psi^\zz(X)$ the function 
$s\mapsto \Tr(Q^{-s}P)$, 
\emph{a priori} defined for $\Re(s)>m+\deg(P)$, extends analytically to $\cz$
with possible simple poles at real integers $k\leq m+\deg P$. The map 
$$P\mapsto \mathrm{Res}_{s=0} \Tr(Q^{-s}P):=\Tr_w(P)$$ 
defines a functional on $\Psi^\zz(X)$
which vanishes on $\Psi^{-\dim(X)-1}(X)$. Moreover, as the notation suggests, 
$\Tr_w$ is independent of the choice of $Q$.
From the definition, it follows easily that $\Tr_w$ vanishes on 
commutators (i.e.\
$\Tr_w$ is a trace), and that it is given by the local expression
\begin{equation}\label{rt}
\Tr_w(P)=\int_{S^*X}p_{[-m]} {\mathcal R}\lrcorner\omega^{m}
\end{equation}
where $p_{[-m]}$ is the homogeneous component of homogeneity $-m$ in the full symbol
of $P$, $S^*X$ is the sphere bundle of $T^*X$, ${\mathcal R}$ is the radial 
vector field on $T^*X$ and $\omega$ is the standard symplectic form. As a consequence, 
the right-hand side of (\ref{rt}) is independent of the quantization defining the full
symbol of $P$. 

\begin{corollary}\label{cme}
The functions $\overline{\zeta}(\delta,s)$, $\oeta(\delta,s)$, 
respectively the family of $1$-forms defined in
\textup{(\ref{fcbf})}, admit analytic extensions to $\cz$ with possible simple poles at 
integers $\dim(N)-2\nz$, respectively $\dim(N)-1-2\nz$.
\end{corollary}
\dem 
The claim follows  directly from Proposition \ref{tilp} with the exception of
possible poles at $\dim N-2\nz-1$, respectively $\dim N-2-2\nz$ and
possible additional poles introduced by the Gamma factors. 
The first type of poles does not arise because
of the (anti)-symmetry of polynomial functions. The poles of $\Gamma(s/2)$ 
occur at $s\in -2\nz$. At these values the
zeta function is finite, since its residue equals the residue trace of a differential
operator and hence vanishes. The function $\Gamma(s+1/2)$ introduces additional 
poles at $s\in -2\nz-1$. The eta function is finite at these points, since its residue
equals $\Tr_w(\delta^{-2k+1})=0$. For the Bismut-Freed connection form, the poles of 
$\Gamma(s/2+1)$ occur at $s=2k-2$ for $k$ a negative integer, and one sees 
that $\Tr((D_U^*D_U)^{-\frac{s}{2}} D_U^{-1}\tna(D_U))$ is also regular at these points.
\qed

\vspace{0.5cm}\section{The determinant line bundle over $S^1$}\label{det}

The determinant line bundle $\det(D)$ is a complex line bundle associated to every
family of elliptic operators. It has a canonical connection due to Bismut and Freed 
\cite{bf2}. In the case where the base is $S^1$ this bundle is trivial, like every
complex vector bundle. Then it is possible to find explicitly 
the connection form in special trivializations. Using (\ref{fhol}) we will get 
a formula for the holonomy along $S^1$ without having to deal with open 
covers; rather, the cover is encoded in the connection form.

Let $\cE^\pm$ be the infinite-dimensional vector bundles over $S^1$
whose fiber over $x\in S^1$ is the space of smooth sections in $E^\pm$ over
$N_x$. For $\alpha\in\rz_+$ let $ V_\alpha:=\{x\in S^1; \alpha\notin 
\Spec(D^*D)\}$. For every $x\in S^1$, let $\cE^\pm_{x}=\cE^\pm_{x,\leq\alpha}
\oplus\cE^\pm_{x,>\alpha}$ be the decomposition of
$\cE^\pm_{x}$ in subspaces spanned by eigensections of $D^*D$, 
respectively of $DD^*$, of eigenvalue smaller or strictly larger than $\alpha$.
Then $\cE^\pm_{<\alpha}$ form vector bundles over $ V_\alpha$ but in general
not over $S^1$.
Nevertheless, if we define
$\det(D)_\alpha\to V_\alpha$ as the line bundle 
$\Lambda^\topp(\cE^+_{<\alpha})^*
\otimes\Lambda^\topp(\cE^-_{<\alpha})$ then $\det(D)_\alpha$ is isomorphic to
$\det(D)_\beta$ over $V_\alpha\cap V_\beta$ via the 
determinant of the linear isomorphism of finite dimensional vector spaces 
$$D:\cE^+_{<\alpha,>\beta}\to \cE^-_{<\alpha,>\beta}.$$
These isomorphisms clearly fulfill the cochain condition so one obtains a line 
bundle over $S^1$ called $\det(D)$. 

$\det(D)_\alpha$ inherits the connection $\tna$ defined in (\ref{cbf}).
The Bismut-Freed connection on $\det(D)_\alpha$ is defined as
$$\nabf_\alpha:=\tna+A_\alpha(0)$$
where $A_\alpha(0)$ is the regularized value in $s=0$ of the
meromorphic extension of the family of $1$-forms
$$A_\alpha(s):=\Gamma\left(1+\frac{s}{2}\right)\Tr_{>\alpha}
((D^*D)^{-\frac{s}{2}} 
D^{-1}\tna(D))$$
which is well defined and analytic for $\Re(s)>n$. One checks easily that
$\nabf_\alpha$ is mapped to $\nabf_\beta$ under the isomorphism 
$\det(D)_\alpha\to\det(D)_\beta$, so we finally get a connection $\nabf$
on $\det(D)$.
 
\begin{lemma} \label{tc}
There exist trivial finite-dimensional vector bundles
$U^\pm$ over $S^1$ and maps of vector bundles $D_{12}:U^+\to \cE^-$, 
$D_{21}:\cE^+\to U^-$, 
$D_{22}:U^+\rightarrow U^-$, such that the operator
\begin{equation}\label{du}
D_U:=\left[\begin{array}{cc} D&D_{12}\\D_{21}&D_{22}\end{array}\right]
:\cin(S^1,\cE^+\oplus U^+)\to\cin(S^1,\cE^-\oplus U^-)
\end{equation}
is invertible; moreover,
there exists a finite cover $\cV=\{V_{\alpha_1}, \ldots,  
V_{\alpha_k}\}$ of $S^1$ such that for 
$x\in V_\alpha\in\cV$, $D_U$ coincides with 
$D$ on ${\mathcal
E}^+_{x,>\alpha}$, and $D_U$ maps $\cE^+_{x,<\alpha}\oplus U^+$
onto $\cE^-_{x,<\alpha}\oplus U^-$. 
\end{lemma}
\dem Let $\{V_{\beta_1}, \ldots, V_{\beta_l}\}$ be an open cover of $S^1$
by sets of the form $V_\alpha$, and $\{V_1', \ldots, V_l'\}$
a compact subcover. Let $\{\phi_j\}$ be a partition of unity
with $\supp(\phi_j)\subset V_{\beta_j}$ and $\phi_j\equiv 1$ on $V_j'$.
Fix isomorphisms $f_j:\cE^+_{<\beta_j}\to \cz^{r_j}$ over $V_{\beta_j}$.
Let $U^-:=\cz^{\sum r_j}$ and $D_{21}:=\sum \phi_j f_j:\cE^+\to U^-$.
By construction, $D\oplus D_{21}$ is injective so its finite dimensional 
cokernel (called $U^+$) forms a vector bundle over $S^1$. Let 
$D_{12}\oplus D_{22}:U^+\to \cE^-\oplus U^-$ be the inclusion map. Then
clearly $D_U$ defined by (\ref{du}) is invertible.
Choose the open cover $\cV$ such that $\max(\beta_j)<\min(\alpha_i)$. 
The second statement of the lemma is easily checked.
\qed

We trivialize $U^+$ and we endow $U^\pm$ with the trivial connections, 
denoted for simplicity $\tna$ (so the connections on $U^+$ and on 
its image in $\cE^-\oplus U^-$ are different). 

\begin{proposition} \label{dlt}
The determinant line bundle with the Bismut-Freed connection is isomorphic to
the trivial bundle $\cz\times S^1\rightarrow S^1$ with the connection
$d+A(0)$, where $A(0)$ is the regularized value in $s=0$ of the meromorphic family of
$1$-forms
\begin{equation}\label{fcbf}
A(D_U,s):=\Gamma\left(1+\frac{s}{2}\right)\Tr((D_U^*D_U)^{-\frac{s}{2}} 
D_U^{-1}\tna(D_U)).
\end{equation}
\end{proposition}
\dem Let us compare the connections $\nabf_{\alpha_j}$ corresponding 
to the operators $D$ and $D_U$, for $\alpha_1,\ldots,\alpha_k$ as in Lemma
\ref{tc}. The meromorphic families of $1$-forms $A_{\alpha_j}(s)$ coincide 
by the second part of Lemma \ref{tc}. The bundles $\det(D)_{\alpha_j}$ and
$\det(D_U)_{\alpha_j}$ are canonically isomorphic together with 
their connections, since $U^\pm$ are trivial with trivial connections.
Thus $(\det(D),\nabf)\simeq(\det(D_U),\nabf)$. Now $\det(D_U)$
can also be constructed using the cover $\{V_0\}$ because $D_U$ is invertible.
For this cover we have $\det(D_U)_0=\cz$ with trivial induced
connection $\tna=d$, and $A_0(s)$ is given by (\ref{fcbf}).
\qed

In the case where the family $D$ is invertible, one obtains a canonical family 
of $1$-forms $A(D,s)$ by taking $U^\pm=0$. 

The Bismut-Freed connection preserves the Quillen metric on $\det(D)$ \cite{quil}. 
In the complex setting, it is just the Hermitian connection associated to that 
metric. We will call its holonomy around $S^1$ the holonomy of $\det(D)$. A simple
computation shows
\begin{equation}\label{fhol}
\hol(\det(D))=\exp\left(-\int_{S^1}A(0)\right).
\end{equation}

The results of this section can be adapted to the case of an arbitrary 
base manifold:
\begin{proposition}
Let $D:E^+\to E^-$ be a family of (classical pseudo-)diffe\-rential elliptic
operators of positive order
on the fibers of a fibration $N\to B$ of compact manifolds.
There exist finite-dimensional vector bundles
$U^\pm\to B$ and maps of vector bundles $D_{12}:U^+\to \cE^-$, 
$D_{21}:\cE^+\to U^-$, 
$D_{22}:U^+\rightarrow U^-$, such that the operator $D_U$ defined in 
\textup{(\ref{du})} is invertible.
There exists a finite cover $\cV=\{V_{\alpha_1}, \ldots,  
V_{\alpha_k}\}$ of $B$ such that for 
$x\in V_\alpha\in\cV$, $D_U$ coincides with 
$D$ on ${\mathcal
E}^+_{x,>\alpha}$, and $D_U$ maps $\cE^+_{x,<\alpha}\oplus U^+$
onto $\cE^-_{x,<\alpha}\oplus U^-$. Endow $U^\pm$ with connections $\tna$, then
$$(\det(D),\nabf)\simeq\left(\cz, d+A(D_U,0)\right)\otimes\left(
\Lambda^\topp(U^+)^*
\otimes\Lambda^\topp (U^-),\tna\right)^{-1}.$$
\end{proposition}
The proof is a simple adaptation of Lemma \ref{tc} and Proposition \ref{dlt}.
We will not use this statement in the present paper.

\vspace{0.5cm}\section{The extended adiabatic algebra}\label{ea}

Adiabatic operators were introduced in \cite{mame90} for the study of 
the adiabatic limit of the Hodge cohomology groups. By a different 
approach using differentiable groupoids a slightly smaller calculus appeared
in \cite[Example 7]{nwx}. We refer to 
\cite{mordis,ha} for a detailed discussion 
of the scalar adiabatic algebra; the extension to bundles is straightforward.

In this section we will introduce 
a larger calculus $\Psi_{ae}$ containing $\pa$ as a subcalculus. This 
calculus realizes the operators of the type we are 
interested in and their complex powers \emph{at the limiting value $t=0$}. 
Therefore, it will suffice to study these limiting operators, which belong to
the so-called \emph{suspended} algebra \cite{meleta}. 

Let $p:X\rightarrow M$ be a fibration of closed manifolds (in this paper, 
we will use $N\rightarrow S^1$), and $E^+, E^-\rightarrow X$ two vector 
bundles. An \emph{adiabatic vector field} is a family $v:[0,\infty)\rightarrow
{\mathcal V}(X)$, such that $v(0)$ is tangent to the fibers of $p$. These 
vector fields form a locally free $C^\infty(X\times[0,\infty))$-module $\acV$,
so they are the sections of a vector bundle over $X\times[0,\infty)$ called $\atx$.
Let $\tilde{E}^\pm$ denote the lifts of $E^\pm$ to $X\times[0,\infty)$ and $U$ the
universal enveloping algebra functor.
An adiabatic family of differential operators from $E^+$ to $E^-$ is an element in 
${\mathrm{Hom}}(\tilde{E}^+,\tilde{E}^-)\otimes_{C^\infty(X\times[0,\infty))}
U(\cin(\atx))$, i.e.\ a composition of adiabatic vector fields and 
bundle homomorphisms.

Let $X^2_a=[X^2\times[0,\infty); X\times_M X\times\{0\}]$ be the blow-up of 
the fiber diagonal of $p$ at $t=0$ inside $X^2\times[0,\infty)$. This means
replacing $X\times_M X\times\{0\}$ by the half-sphere bundle of its positive 
normal bundle, and then gluing along geodesic rays to get the smooth 
structure. The result is a smooth manifold with corners of codimension $2$, 
and the smooth structure is independent of the metric used to 
define geodesics. 

Denote by $\beta$ the canonical blow-down map from
$X^2_a$ to $X^2\times[0,\infty)$.
The \emph{lifted diagonal} $\Delta_a$ is by definition the closure of 
$\beta^{-1}(\Delta_X\times (0,\infty))$ inside $X^2_a$. It is easy to see that
$\beta:\Delta_a\to \Delta_X$ is an diffeomorphism.
Let $p_1,p_2:X^2_a\to X\times [0,\infty)$ be the composition of $\beta$
with the projections on the first, respectively on the second $X$ factor. 

\begin{lemma}[\cite{ha}]\label{lff}
The interior of the 
\emph{front face} $\ffa$ of $X^2_a$ introduced by the blow-up 
is canonically diffeomorphic to $X\times _M X\times_M \atmo$ and therefore has a
canonical vector bundle structure. 
\end{lemma}

One can view adiabatic families of differential operators as distributions 
on $X^2\times[0,\infty)$ having Schwartz 
kernels conormal to $\Delta_X\times(0,\infty)$ and supported on it,
but, because of the degeneracy, not all such kernels define adiabatic 
operators. In fact, the 
Schwartz kernels of adiabatic operators lift to the manifold $X^2_a$ defined 
above, and these lifts span \emph{all} the distributions conormal to and 
supported on $\Delta_a$, and extendable across $\ffa$. 
One must make precise the density bundle where these lifted kernels 
take values; this is $p_2^*\Omega(\atx)$, where 
$\Omega$ is the $1$-density functor. In light of this fact, 
we make the following definition.

\begin{definition}
An adiabatic pseudo-differential operator is a distributional section in
$E^-\boxtimes (\Omega({}^a\!TX)\otimes (E^+)^*)$
over $X^2_a$, classical conormal to $\Delta_a$, extendable across the front face and
vanishing rapidly to all other boundary faces.
\end{definition}

For $s\in\cz$, we denote the space of adiabatic operators of conormality order $s$ 
by $\pa^s(X,E^+,E^-)$.

In Section \ref{gen} we will need a larger calculus that we introduce now. Let 
$U^\pm\rightarrow M$ be vector bundles. Let 
\begin{eqnarray*}
S_{12}&=&[X\times M\times [0,\infty);X\times_M M\times\{0\}]\\
S_{21}&=&[M\times X\times [0,\infty);M\times_M X\times\{0\}]
\end{eqnarray*}
be manifolds with corners obtained through the blow-up of "diagonal copies" 
of $X$. We call $\ff_{21}$, respectively $\ff_{12}$ the new faces introduced by 
blow-up. Let $\Psi^{-\infty}_{12}(X,U^-,E^+)$ be the space of 
smooth sections over $S_{12}$  in the bundle 
$E^+\boxtimes (\Omega({}^a\!TM)\otimes (U^-)^*)$, which are
 rapidly vanishing to the boundary faces other than $\ff_{12}$.
We similarly define $\Psi^{-\infty}_{21}(X,E^-,U^+)$ as the space of smooth 
sections of $p_1^*(U^+)\boxtimes p_2^* (\Omega({}^a\!TX)\otimes (E^-)^*)$ 
over $S_{21}$
which are rapidly vanishing at all boundary faces other than  
$\ff_{21}$. Further, 
define 
$$\Psi_{22}(M,U^+,U^-):=\pa(M,U^+,U^-)$$ 
to be the space of adiabatic 
operators corresponding to the identity fibration $M\rightarrow M$, and 
$\Psi_{11}(X,E^+,E^-):=\pa(X,E^+,E^-)$.

\begin{theorem}\label{compa}
There exists natural composition maps for adiabatic operators:
$$\Psi_{ij}^z\circ\Psi_{jk}^w\subset \Psi_{ik}^{z+w}.$$
\end{theorem}
\dem In this proof we omit the bundles from the notation. 
We will construct appropriate \emph{adiabatic triple 
spaces}. The result will follow from the properties of 
pull-back, product and push-forward operations on conormal distributions
\cite{mel1,mela}. 

Let us first prove the theorem for the adiabatic algebra $\pa(X)$. 
Let $X^3_a$ be the iterated blow-up of
$X^3$,  
$$X^3_a:=[X^3\times[0,\infty); \cF_3\times\{0\}; \cF_{12}\times\{0\},
\cF_{23}\times\{0\}, \cF_{13}\times\{0\}]$$
where 
\begin{eqnarray*}
\cF_3&=&\{(x_1,x_2,x_3)\in X^3; p(x_1)=p(x_2)=p(x_3)\},\\
\cF_{ij}&=&\{(x_1,x_2,x_3)\in X^3; p(x_i)=p(x_j)\},\ \ i,j=1,2,3
\end{eqnarray*}
are fiber diagonals. We claim that there exist p-fibrations
$p_{12},p_{23},p_{13}:X^3_a\to X^2_a$, making the diagrams
\begin{equation}\label{diac}
\xymatrix{X^3_a\ar[d]^{p_{ij}}\ar[r]^-{\beta} &X^3\times[0,\infty)\ar[d]^{p_{ij}}\\
X^2_a\ar[r]^-\beta&X^2\times[0,\infty)}
\end{equation}
commute. By symmetry, it is enough to show this for
$i=1,j=2$. Then $p_{12}$ is defined as the composition of blow-down maps, 
isomorphisms and projections
\begin{eqnarray*}
X^3_a&=&[X^3\times[0,\infty); \cF_3\times\{0\}; \cF_{12}\times\{0\},
\cF_{23}\times\{0\}, \cF_{13}\times\{0\}]\\
&\to& [X^3\times[0,\infty); \cF_3\times\{0\}, \cF_{12}\times\{0\}]\\
&\simeq& [X^2\times X\times[0,\infty); \cF_{12}\times\{0\}; \cF_3\times\{0\}]\\
&\to&[X^2\times X\times[0,\infty); \cF_{12}\times\{0\}]\\
&\simeq&[X^2\times[0,\infty); X\times_M X\times\{0\}]\times X\\
&\to&X^2_a.
\end{eqnarray*}
The commutativity of the blow-up \cite[Prop. 5.8.1]{mela} shows that 
(\ref{diac}) commutes. 

One sees by continuity from the interior that 
\begin{eqnarray}
p_{12}&:&p_{13}^{-1}(\Delta_a)\to X^2_a,\label{uun}\\
p_{13}&:&p_{12}^{-1}(\Delta_a)\cap p_{23}^{-1}(\Delta_a)\to\Delta_a\label{cd}
\end{eqnarray}
are diffeomorphisms. For $A,B\in\pa(X)$ define
\begin{equation}\label{comp}
A\circ B=(p_{13})_*(p_{12}^*A\cdot p_{23}^*B).
\end{equation}
Clearly $p_{12}^*A$ is conormal to $p_{12}^{-1}(\Delta_a)$.
The product is well-defined, since $p_{12}^*A$ and $p_{23}^*B$
are conormal to transverse submanifolds of $X^3_a$. The rapid vanishing 
of adiabatic operators to boundary faces other than $\ffa$
shows that the push-forward is well-defined.
Moreover, from (\ref{uun}), (\ref{cd}) and \cite[Prop. 6.11.5]{mela}, 
the push-forward is again a (classical) conormal distribution to $\Delta_a$.

The composition rule for $\pa(M)$ is a particular case of what we have 
just shown. The remaining six cases are essentially easier, since at least
one of the terms involved is a smooth distribution, but less standard because of 
the non-symmetry. Let us only show, for instance, how to compose
$A\in\pa^z(X)$ with $B\in\Psi_{12}^{-\infty}(X)$. Define
$$S_{112}=[X\times X\times M\times[0,\infty); \cG_3\times\{0\}; 
\cG_{12}\times\{0\}, \cG_{23}\times\{0\}, \cG_{13}\times\{0\}]$$
where 
\begin{eqnarray*}
\cG_3&=&\{(x_1,x_2,x_3)\in X^2\times M; p(x_1)=p(x_2)=x_3\},\\
\cG_{12}&=&\{(x_1,x_2,x_3)\in X^2\times M; p(x_1)=p(x_2)\},\\
\cG_{j3}&=&\{(x_1,x_2,x_3)\in X^2\times M; p(x_j)=x_3\}, j=1,2.
\end{eqnarray*}
There exist p-fibrations $p_{12}:S_{112}\to X^2_a$, $p_{j3}:S_{112}\to S_{12}$,
$j=1,2$, which commute with the blow-down maps 
$S_{112}\to X\times X\times M\times[0,\infty)$, $X^2_a\to X^2\times[0,\infty)$,
$S_{12}\to X\times M\times[0,\infty)$ and the corresponding projections
by analogy with (\ref{diac}). Define the composition of $A,B$ by (\ref{comp}).
There is no issue now with the product of distributions, since $B$ is smooth. The only thing we 
need to note is that $p_{13}$ maps $p_{12}^{-1}(\Delta_a)$ 
diffeomorphically onto $S_{12}$, thus the conormality is "integrated out" through
push-forward and the result is smooth.
\qed

The above theorem is yet another materialization of Melrose's program of 
"microlocalizing boundary fibration structures" \cite{mel1}. Note that 
the composition of operators in $\pa$ appears already in \cite{mame90}, and
for compactly supported adiabatic operators it appears also in \cite{nwx}.
As a consequence, we get a pseudo-differential calculus 
$\Psi_{ae}(X, E^\pm, U^\pm)$ formed by matrices
$$\left[\begin{array}{cc} A_{11}&A_{12}\\A_{21}&A_{22}\end{array}\right],$$
where $A_{ij}$ belongs to $\Psi_{ij}$. 

\begin{remark}\label{psaet}
The algebra $\cin([0,\infty))$ is
central in $\Psi_{ae}$, so we can restrict adiabatic operators to
any given $t$.
For $t>0$, $A_{11}(t)\in\Psi(X,E^+,E^-)$, 
$A_{22}(t)\in\Psi(M,U^+,U^-)$, while the off-diagonal components 
are smoothing operators 
\begin{eqnarray*}
A_{21}(t):C^{-\infty}(M,U^+)&\to&\cin(X,E^-),\\
A_{12}(t):C^{-\infty}(X,E^+)&\to&\cin(M,U^-).
\end{eqnarray*}
Thus for $t>0$, $A(t)$ belongs to an extension of the calculi 
of pseudo-differential operators on $X$
and $M$, independent of $t$, that we call $\Psi_{ex}(X;M)$.
\end{remark}

For $E^+=E^-=E$ and $U^+=U^-=U$, the space ${\cA}=
\Psi_{ae}^\zz(X, E, U)$ is an algebra.
$\cA$ has a natural double filtration 
$\{\cA^{i,k}\}_{i\in\zz,k\in{\mathbb N}}$, increasing in the first 
and decreasing in the second superscript, i.e.\  ${\cA}^{i_1,k_1}
\subset{\cA}^{i_2,k_2}\Leftrightarrow i_1\leq i_2$ and $k_1\geq k_2$. 
Namely, $\cA^{i,k}=t^k \cA^i$, where $i$ denotes the conormality order
while $k$
indicates the order of vanishing at $t=0$. This filtrations are compatible 
with the product on ${\cA}$ by Theorem \ref{compa} (for the first filtration)
and since $\cin([0,\infty))$ is central in $\cA$ (for the second).

\begin{lemma} 
The quotient ${\cA}^{k}/{\cA}^{k-1}$ 
is canonically isomorphic to 
$$S^{[k]}({}^a\!T^* X\setminus 0, E)\oplus 
S^{[k]}({}^a\!T^* M\setminus 0, U),$$
where $S^{[k]}$ denotes symbols of pure homogeneity $k$.
\end{lemma}
\dem Clearly $\Psi_{12}$ and $\Psi_{21}$ die in the quotient.
On $\pa(X)$, the isomorphism is induced by conormal principal symbol short 
exact sequence 
$$0\to\pa^{k-1}(X,E)\hookrightarrow \pa^{k}(X,E)\to 
S^{[k]}(N^*\Delta_a\setminus \{0\},E)\to 0,$$
keeping in mind that the conormal bundle to $\Delta_a$ is $\atsx$. The 
second component comes from $\pa(M)$, which is just a special case of the first.
\qed 

Melrose \cite{meleta} introduced a 
notion of suspension for pseudo-differential operators. If $W$ is any vector bundle 
over $M$, one can define $\Psi_{\sus(W)}(X/M)$ as the space of families of
pseudo-differential operators over the fibers of the fiber bundle $X\times_M
W\rightarrow M$, which are translation invariant with respect to $W$, and whose
Schwartz kernels decay rapidly away from the $W$-diagonal. These operators are
identified (via Fourier transform in $W$) with their \emph{indicial families},
i.e.\ families over $M$ of parameter-dependent operators on the fibers of 
$X\rightarrow M$, with parameters in $W^*$ having symbolic behavior. 

\begin{lemma}
The quotient ${\cA}_\partial:={\cA}^{\zz,0}/{\cA}^{\zz,1}$ 
is canonically isomorphic to the algebra of matrices 
$(a_{ij})_{1\leq i,j\leq2}$ such that
\begin{enumerate}
\item $a_{11}\in \Psi_{\sus({}^a\!T^* M_{|t=0})}(X/M,E)$. The suspending
variables $\tau$ live naturally on ${}^a\!T^* M_{|t=0}$, which is 
canonically isomorphic to $T^* M$. 
\item $a_{22}\in S^\zz({}^a\!T^* M_{|t=0},U)$, 
the space of classical symbols on the vector bundle ${}^a\!T^* M_{|t=0}$, 
with values in $U\otimes U^*$.
\item $a_{12}\in\cS(X\times_M T^* M,E\boxtimes U^*)$, 
$a_{21}\in\cS(T^* M\times_M X,U\boxtimes E^*)$, where $\cS$ denotes smooth 
sections rapidly vanishing at infinity.
\end{enumerate}
\end{lemma}
\dem
Consider the restriction of an adiabatic operator $A$
to $\ffa$. This operation is well-defined, since the conormality locus
$\Delta_a$ of $A$ is transversal to
$\ffa$, and $A$ is extendable across $\ff_a$. 
The interior of $\ffa$ is canonically identified with 
the total space of the vector bundle
$N(\Delta_a)|_{\Delta_a\cap\ffa}$, and $A|_{\ffa}$
vanishes rapidly at infinity. This means exactly that
$A|_{\ffa} \in \Psi_{\sus(\atsM_{|t=0})}(X/M,E)$.
Thus we can Fourier-transform $A|_{\ffa}$ in each fiber, and  
we define the \emph{normal operator} of $A$ by
$$\cN(A):=\cF(A|_{\ffa}).$$
Note that the density factor in $A|_{\ffa}$ is used in the Fourier transform.
 
The second statement is a particular case of the first (note that 
in the case when the fiber is $0$-dimensional, the indicial
operator of suspended operators is just a classical symbol).

For the third statement, note that $\ff_{12}^\circ\simeq X\times_M T M$
and $\ff_{21}^\circ\simeq T M\times_M X$. The isomorphism of the lemma 
is again given by restriction to the front face followed by Fourier transform.
\qed

The normal operator introduced in the above lemma is surjective, 
multiplicative, and plays the role of a 
"principal boundary symbol". Modulo choices, this symbol can be extended
to a full boundary symbol $q$ with values in the formal 
series algebra ${\cA}_\partial[[t]]$.

\begin{lemma}\label{bi}
The \emph{boundary algebra} ${\cA}^{\zz,{\mathbb N}}/
{\cA}^{\zz,\infty}$ is isomorphic as a vector space to ${\cA}_\partial[[t]]$.
\end{lemma}
\dem
The choices involved in the definition of $q$ are
a Riemannian metric on $M$, a connection in the fibration
$X\rightarrow M$ and the restriction of a covariant derivative 
in $E\rightarrow X$ to horizontal vectors. 
Let us start with $\pa(X,E)$. We use an analog of the so-called 
right quantization, i.e.\ a map
\begin{eqnarray*}
X\times_M X\times_M \atm&\rightarrow& X^2\times[0,\infty)\\
(x_m,y_m,v_m(t))&\mapsto& (x_m, \tilde{c}_{v_m(t)}(y_m)(1),t),
\end{eqnarray*}
where $x_m, y_m$ are points in $X$ over $m\in M$, $v$ is an adiabatic vector 
field, $c_v$ is the geodesic starting at $m$
in the direction of $v$ and $\tilde{c}(y_m)$ is 
the horizontal lift of $c$ starting at $y_m$. This map lifts to 
\begin{equation}\label{ima}
X\times_M X\times_M \atm\stackrel{\imath}{\to} X^2_a
\end{equation}
(note that the restriction of $\imath$ to $\{t=0\}$ is 
just the canonical identification $X\times_M X\times_M \atm_{|t=0}\simeq
\ffa^\circ$ from Lemma \ref{lff}).
The map $q$ is defined by pulling back adiabatic operators
to $X\times_M X\times_M\atm$ via $\imath$, taking their Taylor series at $t=0$ 
and then Fourier transforming in 
the fibers of $\atm_{t=0}$. We need the connection in horizontal directions 
in order to trivialize the pull-back of $E$ over each fiber of 
$\atm\rightarrow M$, so that we can compute Fourier transform.

The above definition of $q$ can be specialized to $\pa(M,U)$. 
For $A\in \Psi_{12}$
we proceed similarly: first, lift the map
\begin{eqnarray*}
X\times_M \atm&\rightarrow& X\times M\times[0,\infty)\\
(x_m,v_m(t))&\mapsto& (x_m,{c}_{v_m(t)}(1),t)
\end{eqnarray*}
to a map $\imath:X\times_M \atm\to S_{12}$. Let $q(A)$
be the Fourier transform of the Taylor series of $\imath^*(A)$ at $t=0$. 
Note that each Taylor coefficient of $\imath^*(A)$ is an Euclidean density
on the fibers of $X\times_M \atm\to X$, with values in $E\boxtimes U^*$; the 
density factor is used in the Fourier transform. The definition of $q$ on 
$\Psi_{21}$ is analogous.
\qed

Choose local coordinates $x_j$ in $M$, and let $\tau_j$ be the dual 
coordinates in $T^*M$. Then 
\begin{equation}\label{qd}
q(it\nabla_{\tilde{\partial}_{x_j}})=\tau_j
\end{equation}
is a constant series in $t$, where $\tilde{\partial}_{x_j}$
is the horizontal lift of $\partial_{x_j}$
and $\nabla$ is the connection used in the construction of $q$. If 
$D$ is a family of differential operators on the fibers of $\phi$ viewed as
an adiabatic operator constant in $t$ then 
$$q(D)=D.$$

\begin{lemma} 
The product structure on
${\cA}^{\zz,0}/{\cA}^{\zz,\infty}$ induced by the 
linear isomorphism $q$ is a
deformation with parameter $t$ of the suspended (fiber-wise)
product:
\begin{equation}\label{prod}
A(x,\tau)*B(x,\tau)=AB+it 
\sum_{j=1}^{\dim M}\frac{\partial A}{\partial
\tau_j}\nabla_{\tilde{\partial}_{x_j}}(B)+O(t^2).
\end{equation}
\end{lemma}
\dem The product $*$ is given by a sequence of bi-differential operators
with polynomial coefficients. Thus, these coefficients are determined 
from \emph{differential} adiabatic operators. Both $\partial_{\tau_j}$ 
and $\nabla_{\tilde{\partial}_{x_j}}$ are derivations on $\adel$
so it is enough to prove (\ref{prod}) on a set of generators. Such generators
are for instance families in $t$ of fiberwise differential operators on $X/M$,
and $it\nabla_{\tilde{\partial}_{x_j}}$. For these 
generators (\ref{prod}) is true without the error term.
\qed

Let $Q\in {\cA}^{1,0}$ be an adiabatic operator
such that $Q(t)$ is elliptic, self-adjoint and positive for all $t>0$ and
$\cN(Q)$ is an elliptic self-adjoint positive family of suspended operators.
We call then $Q$ invertible and positive.
 
\begin{lemma}\cite{ha}
Let $Q\in {\cA}^{1,0}$ be an invertible positive adiabatic operator. 
Then there exists the entire family $\{Q^s\}_{s\in\cz}$ of complex 
powers of $Q$ such that $Q^s\in{\cA}^{s,0}$. 
\end{lemma}
\dem
Follow the proof of \cite{buc} with the boundary symbol and the 
conormal principal symbol replacing the principal symbol map.
\qed

The trace function $t\mapsto{\Tr}(A(t))$ of an 
extended adiabatic operator $A \in \Psi_{ae}^{\alpha,0}(X)$ 
is smooth for $t>0$, provided $\Re(\alpha)<-\dim X$. The behavior at $t=0$ 
of this function is controllable, and thus 
provides us with our main tool to study adiabatic limits.

\begin{proposition}
 Let $A$ be an 
extended adiabatic operator in $\Psi_{ae}^{\alpha,0}(X)$ with $\Re(\alpha)<-\dim X$.
Then the function $t^{\dim M}{\Tr}(A(t))$ is smooth on $[0,\infty)$. Moreover, 
\begin{equation}\label{ase}
{\Tr}(A(t))\sim (2\pi t)^{-\dim M} \int_{T^* M} 
\Tr(q(A)) \omega^{\dim M},
\end{equation}
where $\omega$ is the standard symplectic form on $T^* M$, and $\Tr$
in the right-hand side
is the fiberwise trace of the formal series of families of 
suspended operators $q(A)$.
\end{proposition}
\dem For fixed $t>0$, the trace 
of $A(t)$ is the integral of its distributional kernel over the diagonal. 
We cut off this kernel away from the diagonal and from the front face, 
pull it back 
through the map (\ref{ima}) and Fourier transform it in the fibers of 
$\atm$. Recall that for the identity fibration 
$M\rightarrow M$ there exist two canonical bundle maps
$I,J:\atm\rightarrow TM\times[0,\infty)$. Namely, $I$ is the canonical 
inclusion of adiabatic vector fields inside all vector fields, and $J$ is 
the isomorphism given by $v(t)\mapsto v(t)/t$. We denote by the same letters 
the duals maps. Then $J^*(I_* \omega)=\omega/t$. The 
integral of the restriction of the kernel to the diagonal becomes (after 
Fourier transform) the integral on the total space, with a density factor
$(I_*\omega/2\pi)^{\dim M}$. When we go back to $T^* M$ via $J^*$, 
this becomes $(\omega/(2\pi t))^{\dim M}$. Therefore
$${\Tr}(A(t))= (2\pi t)^{-\dim M} \int_{T^* M} 
\Tr(\cF(\imath^*A(t))) \omega^{\dim M}.$$
From Lemma \ref{bi} we know that $\cF(\imath^*A(t))$ is smooth at $t=0$, 
with Taylor series $q(A)$. The integrand is absolutely integrable for 
$\Re(\alpha)<-\dim X$, so the result follows by Lebesgue's dominated 
convergence theorem.
\qed

We need to apply this result to the meromorphic extension of $\Tr(A(t))$. 
For any meromorphic function $f$ and $z\in\cz$, $k\in\zz$
let $c_{k,z}(f)$ be the coefficient of
$(s-z)^k$ in the Laurent expansion of $f$ near $z$.
 
\begin{proposition}\label{til}
 Let $A(s,t)$ be an entire family of operators in ${\cA}$, such that
$A(s)\in {\cA}^{s,0}$. Then 
for each fixed $k\in\zz$ and $z\in\cz$, the function 
$$[0,\infty)\ni t\mapsto t^{\dim M}c_{k,z}(\Tr(A(t)))$$ 
is smooth, and its Taylor expansion at $t=0$ is given by
$$(2\pi)^{-\dim M} c_{k,z}\left(\int_{T^* M}
\Tr(q(A)) \omega^{\dim M}\right).$$
\end{proposition}
\dem It follows from the proof of Lemma \ref{bl} that if $a(s,t)$ is a smooth
family of entire families of symbols, then the Laurent development of 
$\int_V a(s,t)d\mu dg$ around any complex point $s$ depends smoothly on $t$. In
particular, if $a(s,t)\sim \sum_{j=0}^\infty t^j a_j(s)$ is 
the Taylor expansion at $t=0$ of $a(s,t)$ then $c_{k,z}(a(t))\sim\sum_{j=0}^\infty 
t^j c_{k,z}(a_j)$.
We apply this to the full adiabatic symbols $a(s,t)$ of the family $A(s,t)$, 
which consist of a symbol on $^a\! T^*\!X$ and a symbol on $^a\! T^*\!M$. 
The result follows from (\ref{ase}) by noting that 
$a_j(s)$ is the full symbol of $q(A)_{[j]}$, where $_{[j]}$ denotes the 
coefficient of $t^j$ in a formal series.
\qed

\vspace{0.5cm}\section{The adiabatic limit: the invertible case} \label{mate}

Throughout this section we assume that the family $D$ is invertible. With this
assumption we only need to use the adiabatic algebra $\Psi_{a}(N)$ and not the
extended calculus.
As a first application of the formalism from Section \ref{ea}, we
reprove a result from \cite{bc}.
\begin{proposition}\label{inveri}
Assume that the family $D$ is invertible. Then there exists $\epsilon>0$ such 
that the operator $\delta_t$ defined in \textup{(\ref{deltat})} is invertible for all 
$0<t<\epsilon$.
\end{proposition}
\dem In the suspended algebra ${\mathcal
A}_\partial=\Psi_{a}^{\zz,0}(N)/\Psi_{a}^{\zz,1}(N)$, the operator
$$\cN(\delta_t^2)=\left[\begin{array}{cc}\tau^2+D^* D&
0\\0&\tau^2+DD^* \end{array}\right]$$
is invertible \cite{meleta} since $D$ is 
assumed to be invertible, so $\cN(\delta_t)$ is also invertible. 
Moreover, $\sigma_1(\delta_t)\in\cin(\assN,\End(E))$
is invertible. Clearly, $\sigma_1(\cN(\delta_t))=\sigma_1(\delta_t)_{|t=0}$.
Let $B_0\in \pa^{-1,0}(N,E)$ be such that $\cN(B_0)=\cN(\delta_t)^{-1}$
and $\sigma_{-1}(B_0)=\sigma_1(\delta_t)^{-1}$. Thus
$\delta_t B_0=I-R_0$, with $R_0\in\pa^{-1,1}(N)$. 
Let $Q\in\pa^{-1,0}(N,E)$ be an adiabatic operator which realizes the 
asymptotic sum $B_0(I+R_0+R_0^2+\ldots)$. 
Then $R(t):=\delta_t Q-I\in\pa^{-\infty,\infty}$. This ideal
is canonically isomorphic to the algebra of smooth
families (indexed by $[0,\infty)$)
of smoothing operators on $N$, which are rapidly vanishing at $t=0$.
Thus, the norm of $R(t)$ as a bounded operator on $L^2(N)$ tends to $0$, which implies
that $\delta_t Q=I+R(t)$ is invertible for sufficiently small $t$.
\qed

From our choice of adiabatic right quantization, it follows that in the algebra
$\pa^{\zz,0}(N)/\pa^{\zz,\infty}(N)\cong{\cA}_\delta[[t]]$, we have
\begin{eqnarray}
q(\delta_t)&=& \left[
\begin{array}{cc} -\tau& D^*\\ D& \tau \end{array}
\right],\label{qdt}\\ 
q(\delta_t^2)&=& \left[
\begin{array}{cc} D^*D+\tau^2& -it\tna_{\ps}D^*\\ 
it\tna_{{\ps}}D& 
DD^*+\tau^2 \end{array}
\right].\label{qdtd}
\end{eqnarray}

Let $f:[0,\infty)\times\cz\rightarrow\cz$ be a family of meromorphic functions on
$\cz$ indexed by $t$. We say that $f$ is a \emph{smooth family of meromorphic
functions} if for every $z\in\cz$ and every $k\in\zz$, the Laurent coefficient
$t\mapsto c_{k,z}(f(t,\cdot))$ of $(s-z)^k$ in $f(s,t)$
is a smooth function of $t$. 

\begin{theorem}\label{eil}
If the family $D$ is invertible then the family of eta functions 
$\oeta(\delta_t,s)$ -- \emph{a priori}
defined for $t>0$ -- extends to a smooth family 
for $t\in[0,\epsilon)$. Moreover, 
\begin{equation}\label{eal}
\lim_{t\rightarrow 0}\oeta(\delta_t,s) =\frac{1}{i\pi}\int_{S^1} A(D,s),
\end{equation}
where $A(D,s)$ is the Bismut-Freed family of $1$-forms defined by 
\textup{(\ref{fcbf})}. 
\end{theorem}
\dem First, $\delta_t$ is invertible and so (\ref{eta}) becomes
$$\oeta(\delta_t,s)=\pi^{-1/2}\Gamma\left(\frac{1+s}{2}\right)
\Tr((\delta_t^2)^{-\frac{s+1}{2}}\delta_t).$$
From Proposition \ref{til} it follows that $t\oeta(\delta_t,s)$ is a
smooth family of meromorphic functions. Further, by (\ref{qdt}) and (\ref{qdtd})
we have
\begin{eqnarray}
\oeta(\delta_t,s)&\sim_{t\rightarrow 0}
&\frac{\Gamma\left(\frac{1+s}{2}\right)}{\sqrt{\pi}} 
\frac{1}{2\pi t}\int_{T^*S^1} \Tr\left(\left[
\begin{array}{cc} D^*D+\tau^2& -it\tna_{{\ps}}D^*\\
it\tna_{{\ps}}D& 
DD^*+\tau^2 \end{array}\right]^{-\frac{s+1}{2}}\right.\nonumber\\
&&\left.\left[
\begin{array}{cc} -\tau& D^*\\ D& \tau \end{array}
\right]\right)d\tau  d\theta. \label{asee}
\end{eqnarray}
where all products and complex powers are in the sense of the product
(\ref{prod}). For the singular term, i.e.\ the 
coefficient of $t^{-1}$, this implies 
\begin{eqnarray*}
\lim_{t\rightarrow 0} t\oeta(\delta_t,s)&=&\frac{\Gamma\left(\frac{1+s}{2}\right)}
{2\pi^{3/2}} \int_{T^*S^1} \Tr\left(-\tau(D^*D+\tau^2)^{-\frac{s+1}{2}}\right)\\&&
+\Tr\left(\tau(DD^*+\tau^2)^{-\frac{s+1}{2}}\right)d\tau  d\theta
\end{eqnarray*}
and this integral vanishes because it is the integral of an odd function of $\tau$. We compute
now the limit of $\oeta(\delta_t,s)$ as $t$ tends to $0$. Since the matrices of operators
$\left[\begin{array}{cc} D^*D+\tau^2& 0\\0& DD^*+\tau^2 \end{array}\right]$ 
and $\left[\begin{array}{cc} -\tau& D^*\\ D& \tau \end{array} \right]$ 
commute modulo $t$, the constant term in $t$ from (\ref{asee}) becomes
\begin{eqnarray}
\lim_{t\to 0}\oeta(\delta_t,s)&=&
\frac{\Gamma\left(\frac{1+s}{2}\right)}
{2\pi^{3/2}} \int_{T^*S^1} \Tr\left(\left[
\begin{array}{cc} D^*D+\tau^2& 0\\0& DD^*+\tau^2 \end{array}
\right]^{-\frac{s+1}{2}}_{[1]}\right.\nonumber\\
&&\left.\left[
\begin{array}{cc} -\tau& D^*\\ D& \tau \end{array}
\right]\right)d\tau  d\theta\nonumber\\&&
+\frac{\Gamma\left(\frac{1+s}{2}\right)}{2\pi^{3/2}} \int_{T^*S^1} -\frac{s+1}{2}
\Tr\left(\left[\begin{array}{cc} D^*D+\tau^2& 0\\0& DD^*+\tau^2 \end{array}
\right]^{-\frac{s+3}{2}}\right.\nonumber\\&&
\left.\left[\begin{array}{cc}0&-i\tna_{{\ps}}D^*\\i\tna_{{\ps}}D&0\end{array}\right]
\left[\begin{array}{cc} -\tau& D^*\\ D& \tau \end{array}
\right]\right)d\tau  d\theta\nonumber\\
&=:&I(s)+J(s)\label{el}
\end{eqnarray}
where the subscript $_{[1]}$ denotes the coefficient of $t^1$ in a formal series.
Since the product involved in the complex powers is the deformed product
(\ref{prod}), the term $I(s)$ does not vanish directly. In fact, it is 
impossible to compute the integrand before taking the trace. Nevertheless, 
the trace can be computed by assuming formally that the operators 
involved commute. Thus,
$I(s)$ has two contributions coming from the diagonal, namely  
\begin{equation}\label{t1}
\frac{\Gamma\left(\frac{1+s}{2}\right)}
{2\pi^{3/2}}\frac{s+1}{2}\frac{s+3}{2}\int_{T^*S^1} -\tau^2 \Tr\left( 
(D^*D+\tau^2)^{-\frac{s+5}{2}}i\tna_{{\ps}}(D^*D)\right)d\tau  
d\theta,
\end{equation}
respectively 
\begin{equation}\label{t2}
\frac{\Gamma\left(\frac{1+s}{2}\right)}
{2\pi^{3/2}}\frac{s+1}{2}\frac{s+3}{2}\int_{T^*S^1} \tau^2 \Tr\left( 
(DD^*+\tau^2)^{-\frac{s+5}{2}}i\tna_{{\ps}}(DD^*)\right)d\tau 
 d\theta.
\end{equation}
These factors cancel each other before integration; alternately, each of them 
equals $0$ after integration; for instance, the quantity (\ref{t1}) equals
$$\frac{\Gamma\left(\frac{1+s}{2}\right)}
{2\pi^{3/2}}\frac{s+1}{2}\int_{T^*S^1} i\tau^2 {\ps}\Tr\left( 
(D^*D+\tau^2)^{-\frac{s+3}{2}}\right)d\tau  
d\theta=0$$
so $I(s)=0$ for large $\Re(s)$. By unique continuation,
\begin{equation}\label{iez}
I(s)\equiv 0.
\end{equation}
We focus now on the second term $J(s)$. It is given by
\begin{eqnarray}
J(s)&=&
-\frac{s+1}{2}\frac{\Gamma\left(\frac{1+s}{2}\right)}{2\pi^{3/2}} \int_{T^*S^1} 
\left[\Tr\left((D^*D+\tau^2)^{-\frac{s+3}{2}}(-i\tna_{{\ps}}D^*)D
\right)\right.\nonumber\\&&
\left.+\Tr\left((DD^*+\tau^2)^{-\frac{s+3}{2}}(i\tna_{{\ps}}D)D^*
\right)\right]d\tau d\theta.\label{inteq}
\end{eqnarray}
Notice that
\begin{eqnarray*}
\lefteqn{-\frac{s+1}{2}
\left[\Tr\left((D^*D+\tau^2)^{-\frac{s+3}{2}}(\tna_{{\ps}}D^*)D
\right)\right.}&&\\&&
+\left.\Tr\left((DD^*+\tau^2)^{-\frac{s+3}{2}}(\tna_{{\ps}}D)D^*
\right)\right]\\
&=&\Tr\left(\tna_{{\ps}}(D^*D+\tau^2)^{-\frac{s+1}{2}}\right)\\
&=&\ps\Tr\left((D^*D+\tau^2)^{-\frac{s+1}{2}}\right)
\end{eqnarray*}
so the \emph{sum} of the two terms in (\ref{inteq}) vanishes. Thus
\begin{eqnarray}
J(s)&=&
-i\frac{\Gamma\left(\frac{3+s}{2}\right)}{\pi^{3/2}} \int_{T^*S^1} 
\Tr\left((DD^*+\tau^2)^{-\frac{s+3}{2}}(\tna_{{\ps}}D)D^*
\right)d\tau d\theta\label{js}\\
&=&-i\frac{\Gamma\left(\frac{3+s}{2}\right)}{\pi^{3/2}}
f\left(\frac{s+3}{2}\right)\int_{S^1} \Tr\left((DD^*)^{-\frac{s}{2}-1}
(\tna_{{\ps}}D)D^*\right)d\theta\nonumber\\&=&
\frac{\Gamma\left(\frac{s}{2}+1\right)}{i\pi}\int_{S^1}
\Tr\left((D^*D)^{-\frac{s}{2}}D^{-1}\tna_{{\ps}}(D)\right)d\theta \label{ja}
\end{eqnarray}
where $f(s)$ is the function from (\ref{fs}). We have used the commutation
formula $(D^*D)^sD^*=D^*(DD^*)^s$.
In conclusion, (\ref{iez}) and (\ref{inteq}) imply (\ref{eal}).
\qed

This result (for $s=0$) is used in the index formula of \cite{phind}
to relate the boundary term with the adiabatic limit of the eta invariant.

\begin{corollary}
The limit of  $\exp(-i\pi\oeta(\delta_t))$ exists and equals the holonomy of $\det(D)$.
\end{corollary}
\dem By specializing to the constant coefficient $\oeta(\delta_t)$
of $\oeta(\delta_t,s)$ at $s=0$, we get
$$\lim_{t\rightarrow 0}e^{-i\pi\oeta(\delta_t)}=e^{-\int_{S^1} A(D,0)}.$$
The right-hand side equals the holonomy along the circle of the connection 
$d+A(D,0)$ in the trivial bundle.
\qed

We turn now to the zeta function. Recall that if the family $D$ is invertible 
then the operator $\delta_t$ defined in (\ref{deltat}) is also invertible for 
small enough $t>0$ (Proposition \ref{inveri}).
\begin{theorem}\label{mz}
Assume that the family of operators $D$ is invertible. Then the family of meromorphic
functions
$t\overline{\zeta}(\delta_t,s)$ extends smoothly to $t=0$. Moreover, 
\begin{equation}\label{lz}
\lim_{t\rightarrow 0}t\overline{\zeta}(\delta_t,s)=\frac{1}{\sqrt\pi}
\int_{S^1}\overline{\zeta}(D,s-1)d\theta,
\end{equation}
while 
$$\frac{d(t\overline{\zeta}(\delta_t,s))}{dt}|_{t=0}=0.$$
\end{theorem}
\dem As for the eta function, by Proposition \ref{til}, (\ref{qdt}) and 
(\ref{qdtd}) we have
\begin{eqnarray}
t\overline{\zeta}(\delta_t,s)&=&t\Gamma\left(\frac{s}{2}\right)
\Tr\left((\delta_t^2)^{-\frac{s}{2}}\right)\nonumber\\
&\sim_{t\to 0}& \Gamma\left(\frac{s}{2}\right)\frac{1}{2\pi}
\int_{T^*S^1} \Tr\left[
\begin{array}{cc} D^*D+\tau^2& -it\tna_{\tps}D^*\\ 
it\tna_{\tps}D& 
DD^*+\tau^2 \end{array}\right]^{-\frac{s}{2}}d\tau d\theta.\label{limz}
\end{eqnarray}
The leading term in this expression is
\begin{eqnarray*}
\lefteqn{\pi^{-1}\Gamma\left(\frac{s}{2}\right)\int_{T^*S^1}\Tr(D^*D+\tau^2)^{-s/2}
d\tau  d\theta}&&\\
&=&\pi^{-1}\Gamma\left(\frac{s}{2}\right)f\left(\frac{s}{2}\right)
\int_{S^1}\Tr(D^*D)^{-\frac{s-1}{2}}d\theta.
\end{eqnarray*}
Together with (\ref{fs}), this implies (\ref{lz}). It is clear that the 
off-diagonal terms in $q(\delta_t^2)$ contribute to (\ref{limz}) only 
modulo $t^2$. It follows that
the coefficient of $t$ in (\ref{limz}) is a sum of two terms coming from the
diagonal. These terms are similar to (\ref{t1})
and (\ref{t2}), only that they are odd in $\tau$ and thus both equal $0$ even before
integrating over $S^1$.
\qed

We can deduce from this result the adiabatic limit behavior of the usual quantities
$\zeta(\delta_t,0)$ and of $\det(\delta_t):=e^{-{\zeta}'(\delta_t,0)}$.
Let $\Tr_\zeta|D|$ denote the regularized value $\zeta(D,-1)$.
\begin{corollary}\label{pzl}
The functions $t\zeta(\delta_t,0)$ and $t\log\det(\delta_t)$ extend smoothly to
$t=0$. Moreover,
\begin{eqnarray}\label{z0}
\lim_{t\rightarrow
0}t\zeta(\delta_t,0)&=&-\int_{S^1}\Tr_w|D|d\theta\\
\lim_{t\rightarrow 0}t\log\det(\delta_t)&=&\int_{S^1}\left(-\Tr_\zeta|D|
+(2\log 2-2)\Tr_w|D|\right)d\theta.\label{ldet}
\end{eqnarray}
\end{corollary}
\dem In terms of $\zeta$, formula (\ref{lz}) reads
\begin{equation}\label{rz}
\lim_{t\rightarrow0}t\Gamma\left(\frac{s}{2}\right)\zeta(\delta_t,s)
=\pi^{-1/2}\Gamma\left(\frac{s-1}{2}\right)\int_{S^1}\zeta(D,s-1)d\theta.
\end{equation}
Recall that
$$\Gamma(0)=-\gamma,\ \ \Gamma'\left(-\frac{1}{2}\right)=2\pi^{1/2}
(\gamma+2\log 2-2),$$
where $\gamma$ is the Euler constant, and $\Gamma(0)$ is the finite part of
$\Gamma(s)$ at $s=0$. Thus
\begin{eqnarray*}
\lefteqn{\lim_{t\rightarrow0}t\left(\frac{2}{s}-\gamma\right)\left(\zeta(\delta_t,0)+
s(\zeta)'(\delta_t,0)\right)}&&\\
&\equiv&
\frac{1}{\sqrt\pi}\left(-2\pi^{1/2}+\pi^{1/2}(\gamma+2\log 2-2)s\right)\\&&
\left(\frac{1}{s}\int_{S^1}\Tr_w(D^*D)^{1/2}d\theta+\int_{S^1}\zeta(D,-1)d\theta\right)
\pmod{s}.
\end{eqnarray*}
Identifying the coefficients of $s^{-1}$ and $s^0$ we obtain (\ref{z0})
and (\ref{ldet}).
\qed

An interesting formula appears in odd dimensions:
\begin{corollary}
Assume that $N$ is odd-dimensional. Then
$$\lim_{t\rightarrow 0}\det(\delta_t)^t
=e^{-\int_{S^1}\Tr_\zeta(|D|)d\theta}.$$
\end{corollary}
\dem
In this case $\zeta(\delta_t,0)$ 
must vanish and by (\ref{z0}) this implies that
$\int_{S^1}\Tr_w(D^*D)^\frac{1}{2}d\theta=0$. The result follows from Corollary
\ref{pzl}.

\vspace{0.5cm}\section{The adiabatic limit: the non-invertible case}\label{gen}

Consider an elliptic family $D$ as in Section \ref{mate}, without the 
invertibility hypothesis. Then $\delta_t$ is also non-invertible in general.
One could define the complex powers of $\delta_t^2$ as being $0$ on the null-space
of $\delta_t$. This works fine for a fixed value of $t$. However, the function 
$(0,\infty)\times\cz\ni(\lambda,s)\mapsto
\lambda^s$ does not converge to $0$ as $\lambda$ tends to $0$, and this means that
the family $(\delta_t^2)^s$ defined in this way is discontinuous at values of $t$ where an eigenvalue
crosses $0$. This turns out to be a major difficulty when trying to extend the
analysis from the previous section to the family $\delta_t$. The
extended adiabatic algebra proves to be the right tool to overcome this problem.

Let $U^+, U^-$ be trivial finite-dimensional 
complex vector bundles over $S^1$, such that there exists an 
invertible family of operators $D_U$ extending $D$ by the formula (\ref{du}). 
As in Section \ref{det}, we extend $\tna_{\ps}$ trivially on
$\cin(S^1,U^\pm)$.
Let $\phi$ be a Schwartz function on 
$\rz$ such that $\phi(0)=1$. Then 
\begin{equation}\label{q}
Q:=\left[\begin{array}{cc} D&D_{12}\phi(\tau)\\ 
D_{21}\phi(\tau)&D_{22}\phi(\tau)\end{array}\right]
\end{equation}
belongs to
$\Psi_{ae}^{1,0}(N;E^+\oplus U^+,E^-\oplus U^-)/
\Psi_{ae}^{1,1}(N;E^+\oplus U^+,E^-\oplus U^-)$. 
There exists $R_t\in\Psi_{ae}^{-\infty}(N)$ such that $q(R_t)=
\left[\begin{array}{cc} 0&D_{12}\phi(\tau)\\ 
D_{21}\phi(\tau)&D_{22}\phi(\tau)\end{array}\right]$
is constant in $t$. By abuse of notation we write $D$ instead of 
$\left[\begin{array}{cc}D&0\\0&0\end{array}\right]$.
Then $q(R_t+D)=Q$. 
Let $E:=E^+\oplus E^-$, $U:=U^+\oplus U^-$ and 
\begin{equation}\label{dt}
d_t=\left[\begin{array}{cc}
-ti\tna_{\ps}&D^*+R_t^*\\D+R_t&ti\tna_{\ps}
\end{array}\right]\in\Psi_{ae}^{1,0}(N, E,U).
\end{equation}
If we replace 
$D$ by $Q$, the identities (\ref{qdt}), (\ref{qdtd}) become
\begin{eqnarray}
q(d_t)&=& \left[
\begin{array}{cc} -\tau& Q^*\\ Q& \tau \end{array}
\right],\label{qdtq}\\ 
q(d_t^2)&=& \left[
\begin{array}{cc} Q^*Q+\tau^2+it\partial_\tau(Q^*)\tna_{\ps}(Q)
& -it\tna_{\ps}Q^*\\ 
it\tna_{{\ps}}Q& 
QQ^*+\tau^2 +it\partial_\tau(Q)\tna_{\ps}(Q^*)\end{array}
\right].\label{qdtdq}
\end{eqnarray}

In particular, the boundary symbol 
$\cN(d_t^2)=\left[\begin{array}{cc}Q^* Q+\tau^2&0\\
0&QQ^*+\tau^2\end{array}\right]$ is invertible inside 
${\cA}_\partial(N,E,U)$ since $Q(0)$ is invertible. 
Proposition \ref{inveri} shows that $d_t$ is 
invertible for small enough $t$. Formulas (\ref{el}), (\ref{js}) from
Section \ref{mate} have analogs with $\delta_t$ replaced
with $d_t$ and $D$ with $Q$. Such results are, of course, 
uninteresting, since $d_t$ and $Q$ depend on the choice of $U^\pm$ and of the 
Schwartz function $\phi$, and moreover $I(s)$, $J(s)$ cannot be simplified 
any further. However, there is 
one instance when these results can be linked to $\delta_t$ and $D$, namely when 
we restrict our attention to the constant coefficient in $s$ at $s=0$. 
The reason is that the values at $s=0$
of the eta and zeta functions of the operators $\delta_t$ and $d_t$ are
closely related.

\begin{lemma} \label{hindep}
Let $d$ be a symmetric elliptic operator of positive order in the algebra 
$\Psi_{ex}(X;M)$ defined in Remark \textup{\ref{psaet}}. Let 
$\mu\mapsto R_\mu\in\Psi_{ex}^{-\infty}(X;M)$ be a $1$-parameter 
family of symmetric smoothing operators. 
Then $\oeta(d+R_\mu)$ is independent of $\mu$
modulo $2\zz$.
\end{lemma}
\dem 
We first show that $\oeta(d+R_\mu)$ is constant in $\mu$ as 
long as $d+R_\mu$ remains invertible. Write $d_\mu:=d+R_\mu$.
Then
\begin{eqnarray*}
\partial_\mu \Tr\left((d_\mu^2)^{-\frac{s+1}{2}}d_\mu\right)&=& 
\Tr\left((d_\mu^2)^{-\frac{s+3}{2}}\left(-\frac{s+1}{2}\right)
(\partial_\mu(R_\mu)d_\mu+d_\mu\partial_\mu(R_\mu))
d_\mu \right.\\ &&
\left.+(d_\mu^2)^{-\frac{s+1}{2}}\partial_\mu(R_\mu)\right)\\
&=&-s\Tr\left((d_\mu^2)^{-\frac{s+1}{2}}\partial_\mu(R_\mu) \right).
\end{eqnarray*}
The last expression is the trace of an entire family of smoothing operators 
so it is entire in $s$,
and it vanishes at $s=0$ because of the $s$ factor. 

Let $\mu_0$ be a point
where $0$ is an eigenvalue of $d_{\mu_0}$. Choose $\alpha\in\rz$ with
$\pm\alpha\notin\Spec(d_{\mu_0})$. Then for $\mu$ close to $\mu_0$,
$\pm \alpha\notin\Spec(d_{\mu})$. Moreover, if $P_\alpha(\mu)$ 
is the projection on the (finite-dimensional) span of eigenspaces of $d_\mu$
of eigenvalue between $\pm\alpha$,
then the map $\mu\mapsto P_\alpha(\mu)\in\Psi_{ex}^{-\infty}(X;M)$ 
is smooth in a neighborhood of $\mu_0$. Clearly, in this neighborhood
$\oeta(d_\mu)-\oeta(d_\mu+P_\alpha)\in 2\zz$, and 
$\oeta(d_\mu+P_\alpha(\mu))$ is constant in $\mu$ by what we 
proved above since $d_\mu+P_\alpha(\mu)$ is invertible.
\qed

\begin{lemma}\label{eeez}
Let $d_t$ be the extended adiabatic operator defined in \textup{(\ref{dt})}
and $\delta_t$ the family of differential operators defined in  \textup{(\ref{deltat})}. 
For all $t>0$, 
\begin{eqnarray*}
\oeta(d_t)&\equiv&\oeta({\delta_t}) +\ind(D)\pmod{2\zz}\\
\zeta(d_t)&=&\zeta(\delta_t).
\end{eqnarray*}
\end{lemma}
\dem
For a fixed $t>0$ consider the operator
$$q_t=\left[\begin{array}{ccc}
\delta_t&0&0\\
0&-it\ps&0\\
0&0&it\ps\end{array}\right]
\in\Psi_{ex}(X;M;E,U)$$
acting on $\cin(N, E^+\oplus E^-)\oplus\cin(S^1,U^+)\oplus\cin(S^1,U^-).$
Note that $q_t, d_t$ act on the same space but are written in basis 
which differ by a permutation. Moreover $d_t-q_t\in\Psi_{ex}^{-\infty}(X;M)$
so by Lemma \ref{hindep}, $\oeta(q_t)-\oeta({d_t})\in 2\zz$. On the other hand,
$\oeta(q_t,s)=\oeta(\delta_t,s)+(\dim U^--\dim
U^+)\oeta(it\ps,s)$.
The result follows by noting that $\oeta(it\ps,s)=1$
 and $\dim U^--\dim U^+=\ind(D)$.

Similarly, for the zeta function 
we see in the notation of Lemma \ref{hindep} that
\begin{eqnarray*}
\partial_\mu \Tr\left((d_\mu^2)^{-\frac{s}{2}}\right)&=&
-\frac{s}{2}\Tr\left((d_\mu^2)^{-\frac{s}{2}-1}
(\partial_\mu(R_\mu)d_\mu+d_\mu\partial_\mu(R_\mu))\right)\\ &=&
-s \Tr\left((d_\mu^2)^{-\frac{s}{2}-1} d_\mu\partial_\mu(R_\mu)\right)
\end{eqnarray*}
equals $s$ times an entire function, so it vanishes at $s=0$. Secondly, there are no
jumps in $\zeta(d_\mu)$ when eigenvalues of $d_\mu$ cross $0$ 
(on finite parts of the
spectrum, $\zeta$ at $s=0$ is just the number of eigenvalues). Thirdly, the function
$\zeta(it\ps,s)$ vanishes at 
$s=0$. Indeed, this is equivalent to the fact that
for the Riemann zeta function,
$$\zeta(0)=-\frac12.$$
\qed

We can now prove our main result for non-invertible families. 
The first formula was proved in \cite{bf2} in the particular case 
when $D$ is a family of compatible Dirac operators. In the second 
formula $|D|$ is defined only on the orthogonal 
complement of the null-space of $D$. 
\begin{theorem}
Let $D$ be a family of elliptic differential operators of order $1$ over $S^1$, and
$\delta_t$ the operator defined by \textup{(\ref{deltat})}. Then 
$\exp(-i\pi\oeta(\delta_t))$ and 
$t\zeta(\delta_t)$ extend smoothly to $t=0$. Moreover,
\begin{eqnarray}
\lim_{t\rightarrow 0} e^{-i\pi\oeta(\delta_t)}&=&(-1)^{\ind(D)}\hol(\det(D))
\label{aleg}\\
\lim_{t\rightarrow 0}
t\zeta(\delta_t,0)&=&-\int_{S^1}\Tr_w|D|d\theta\label{alzg}
\end{eqnarray}
and
\begin{eqnarray}
\lim_{t\rightarrow
0}\left(\zeta(\delta_t,0)+\frac{1}{t}\int_{S^1}\Tr_w|D|d\theta\right)&=&0.
\label{stz}
\end{eqnarray}
\end{theorem}
\dem Working in the extended adiabatic algebra, we will 
see as in Section \ref{mate}
that the families of meromorphic functions
$\oeta(d_t,s)$ and $t\zeta(d_t,s)$ are smooth down to $t=0$. Our strategy 
is to compute the asymptotics at $t=0$ of the regularized value
at $s=0$ of these functions, and then to
use Lemma \ref{eeez} to deduce the asymptotics of $\exp(-i\pi\oeta(\delta_t))$ 
and $t\zeta(\delta_t,0)$.

\paragraph{The eta invariant.}
From Proposition \ref{til} and the definition 
of the eta function we write
\begin{eqnarray}
\oeta(d_t,s)&\sim_{t\rightarrow 0}
&\frac{\Gamma\left(\frac{1+s}{2}\right)}{\sqrt{\pi}} 
\frac{1}{2\pi t}\int_{T^*S^1} \Tr\left(q(d_t^2)^{-\frac{s+1}{2}}
*q(d_t)\right)d\tau d\theta.
\end{eqnarray}
where the products and complex powers are in the sense of the product
(\ref{prod}). Using (\ref{qdtq}), (\ref{qdtdq}) we get
for the coefficient of $t^{-1}$
\begin{eqnarray*}
\lim_{t\rightarrow 0} t\oeta(d_t,s)&=&\frac{\Gamma\left(\frac{1+s}{2}\right)}
{2\pi^{3/2}} \int_{T^*S^1} \Tr\left(-\tau(Q^*Q+\tau^2)^{-\frac{s+1}{2}}\right)\\&&
+\Tr\left(\tau(QQ^*+\tau^2)^{-\frac{s+1}{2}}\right)d\tau  d\theta.
\end{eqnarray*}
Decomposing the trace over a basis of eigensections of $Q^*Q$, $QQ^*$ we 
see that the above terms cancel each other (both can be made zero by choosing 
$\phi$ to be even in $\tau$). By analogy with (\ref{el}) we write

\begin{eqnarray}
\lim_{t\to 0}\oeta(d_t,s)&=&
\frac{\Gamma\left(\frac{1+s}{2}\right)}
{2\pi^{3/2}} \int_{T^*S^1} \Tr\left(\left(
q(d_t^2)^{-\frac{s+1}{2}}\right)_{[1]}\left[
\begin{array}{cc} -\tau& 0\\ 0& \tau \end{array}
\right]\right)d\tau  d\theta\nonumber\\&&
+\frac{\Gamma\left(\frac{1+s}{2}\right)}{2\pi^{3/2}} \int_{T^*S^1} -\frac{s+1}{2}
\Tr\left(\left[\begin{array}{cc} Q^*Q+\tau^2& 0\\0& QQ^*+\tau^2 \end{array}
\right]^{-\frac{s+3}{2}}\right.\nonumber\\&&
\left.\left[\begin{array}{cc}0&-i\tna_{{\ps}}Q^*\\i\tna_{{\ps}}Q&0\end{array}\right]
\left[\begin{array}{cc} -\tau& Q^*\\ Q& \tau \end{array}
\right]\right)d\tau  d\theta\nonumber\\
&=:&I(s)+J(s).\label{elq}
\end{eqnarray}
By integration by parts with respect to $\theta$ we get (see (\ref{js}))
\begin{equation}
J(s)=-i\frac{\Gamma\left(\frac{3+s}{2}\right)}{\pi^{3/2}} 
\int_{T^*S^1} \Tr\left((QQ^*+\tau^2)^{-\frac{s+3}{2}}(\tna_{{\ps}}Q)Q^*
\right)d\tau d\theta.\label{jsq}
\end{equation}
However, $Q$ 
depends on $\tau$ since it involves the function $\phi(\tau)$, so 
the analog of (\ref{ja}) and the explicit 
decomposition of $I(s)$ in (\ref{t1}), (\ref{t2}) fail. 

Let $\cV=\{V_{\alpha_1}, \ldots, V_{\alpha_k}\}$ be the cover of $S^1$ 
from Lemma \ref{tc}, and $\{W_{\alpha_j}\}$ a partition of $S^1$ in 
intervals such that $W_{\alpha_j}\subset V_{\alpha_j}$. We split the integral 
over $S^1$ from the definition (\ref{elq}) of $I(s)$, respectively from 
the formula (\ref{jsq}) for $J(s)$, in the sum of the integrals over 
$W_{\alpha_j}$:
$$I(s)=\sum_{j=1}^k I_{\alpha_j}(s), \ \ J(s)=\sum_{j=1}^k J_{\alpha_j}(s)$$
(it is worth mentioning that integration by parts in $\theta$ was done 
\emph{before} this splitting).
Notice that over $W_{\alpha_j}$ the operators $Q(\tau)$, $\tau$ and
$\partial_\tau Q$ preserve the decomposition 
$$\cE\oplus U=\cE_{>\alpha_j}
\oplus\left(\cE_{<\alpha_j}\oplus U\right)$$
while $\ps Q$ does not necessarily do so. Nevertheless, the form
(\ref{prod}) of the product shows that over $W_{\alpha_j}$ we can compute
$I_{\alpha_j}$, respectively $J_{\alpha_j}$, by projecting $Q$ onto
and taking 
the trace on $\cE_{>\alpha_j}$, resp.\ $\cE_{<\alpha_j}\oplus U$:
\begin{equation}\label{dec}
I_{\alpha_j}(s)=I_{\alpha_j}^>(s)+I_{\alpha_j}^<(s),\ \ 
J_{\alpha_j}(s)=J_{\alpha_j}^>(s)+J_{\alpha_j}^<(s).
\end{equation}

One virtue of this decomposition is that $Q_{>\alpha_j}$, and so also
$I_{\alpha_j}^>(s)$ and $J_{\alpha_j}^>(s)$, do not depend on $\phi$ anymore. 
Thus as in Theorem \ref{eil}, $I_{\alpha_j}^>(s)$ has two components (see
(\ref{t1}), (\ref{t2})) which cancel each other (they are also exact 
forms in $\theta$ but we cannot deduce that their integral vanishes
individually), while for $J_{\alpha_j}^>(s)$, (\ref{jsq})
implies the analog of (\ref{ja})
$$J_{\alpha_j}^>(s)=\frac{1}{i\pi}\int_{W_{\alpha_j}} A(D,s)_{>\alpha_j}.$$

A second feature of (\ref{dec}) is that $I_{\alpha_j}^<(s)+J_{\alpha_j}^<(s)$
involves only traces of finite-dimensional linear endomorphisms. More precisely,
$I_{\alpha_j}^<(s)$, $J_{\alpha_j}^<(s)$ are given by the first term
in (\ref{el}), respectively by (\ref{js}), where the integral is computed 
over $W_{\alpha_j}\times\rz$,  $D$ is replaced by $Q_{<\alpha_j}$ (here we note
that $Q_{<\alpha_j}(\theta, \tau)$ is a classical symbol on 
$W_{\alpha_j}\times\rz$, invertible for $\tau=0$), 
and hence the trace $\Tr$ becomes the endomorphism trace $\tr$. 
Let us examine what happens
when we modify the function $\phi$, of course with the restriction that
$\phi(0)>0$. Lemma \ref{hindep} shows that $\oeta(d_t)$ (hence also 
its limit as $t\to 0$) is unaffected. We have seen that
$I_{\alpha_j}^>(s)$ and $J_{\alpha_j}^>(s)$ do not involve $\phi$. 
We deduce that
$$L^<(s):=\sum_{j=1}^k I_{\alpha_j}^<(s)+J_{\alpha_j}^<(s)$$
is also independent of $\phi$ at $s=0$, 
as long as $\phi$ is a Schwartz function.
However, $L^<(s)$ makes sense and is meromorphic in $s$
for any $\phi$ a classical symbol of 
order $0$ on $\rz$. Indeed, for such $\phi$ all operators involved in 
$I_{\alpha_j}^<(s)+J_{\alpha_j}^<(s)$ are still families of 
classical symbols on $\rz$ indexed by $W_{\alpha_j}$, 
so Lemma \ref{bl} applies.

Let us compute $L^<(s)$ with $\phi$ replaced by $1$. 
Then the corresponding $Q_{\alpha_j}$
are $\tau$-free. The same argument as above for $I_{\alpha_j}^>(s)$ and 
$J_{\alpha_j}^>(s)$ (which was explained in Theorem \ref{eil})
shows that 
$$
I_{\alpha_j}^<(s)=0,\ \ J_{\alpha_j}^<(s)=\frac{1}{i\pi}
\int_{W_{\alpha_j}} A(D,s)_{<\alpha_j}.$$ 
We justify below this substitution (this completes the proof of (\ref{aleg})).

Both $I_{\alpha_j}^<(s)$ and $J_{\alpha_j}^<(s)$ are continuous with 
respect to variations of the Schwartz function
$\phi$. We define a deformation for $0< \mu\leq 1$ by
$\phi_\mu(\tau):=\phi(\tau \mu)$. Note that $\phi_1=\phi$ and $\phi_\mu\to 1$ 
pointwise as $\mu\to 0$. Moreover, 
$|\phi_\mu(\tau)|\leq 1$  
while $|\partial_\tau\phi_\mu(\tau)|<C(1+|\tau|)^{-1}$, uniformly in $\mu$.
We aim to find uniform $L^1$ upper bounds for the integrands in 
$I_{\alpha_j}^<(s)(\mu)$ and $J_{\alpha_j}^<(s)(\mu)$.
First, we note
\begin{eqnarray*}
\|Q_{<\alpha_j}\|&\leq&C\\
\|\ps (Q_{<\alpha_j})\|&\leq&C\\
\|\partial_\tau (Q_{<\alpha_j})\|&\leq&C(1+|\tau|)^{-1}\\
\|({Q_{<\alpha_j}}{Q_{<\alpha_j}^*}+\tau^2)^{-s}\|&\leq&
(\psi(\tau)+\tau^2)^{-\Re(s)}
\end{eqnarray*}
where the last inequality holds for $\Re(s)\geq 0$
with some compactly supported non-negative function $\psi$ with $\psi(0)>0$,
independent of $\mu$.
We will use the inequality 
$$|\tr(A)|\leq l\|A\|$$
for an endomorphism $A$ of an $l$-dimensional vector space. The integrand in
$J_{\alpha_j}^<(s)(\mu)$ was defined using (\ref{jsq}). Using the above bounds
we get
\begin{eqnarray*}
|\tr\left(({Q_{<\alpha_j}}{Q_{<\alpha_j}^*}+\tau^2)^{-\frac{s+3}{2}}
(\tna_{{\ps}}{Q_{<\alpha_j}}){Q_{<\alpha_j}^*}\right)|
& \leq& C_j(\psi(\tau)+\tau^2)^{-\frac{\Re(s)+3}{2}}
\end{eqnarray*}
from which we retain that the integrand in $J_{\alpha_j}^<(s)(\mu)$
is bounded (uniformly in $\mu$) 
by an $L^1$ function of $\tau$ for all $s$ with $\Re(s)>-2$,
in particular for $s=0$.

The integrand in $I_{\alpha_j}^<(s)(\mu)$ admits a similar bound. We show this
for the term coming from the upper left corner, the other one 
being entirely similar. Using (\ref{qdtdq}) we first isolate
the term coming from $it\partial_\tau(Q_{<\alpha_j})\tna_{\ps}(Q_{<\alpha_j}^*)$.
This term is bounded uniformly in $\tau$ for $\Re(s)>-3$ by 
$$C_j(s)(\psi(\tau)+\tau^2)^{-\frac{\Re(s)+3}{2}}(1+|\tau|)^{-1}\tau,$$
which is $L^1$ for $\Re(s)>-2$ as before. Finally, the remaining term
coming from
$\left((Q_{<\alpha_j}
Q_{<\alpha_j}^*+\tau^2)^{-\frac{s+1}{2}}\right)_{[1]}\tau$,
though inexplicit, is bounded for $\Re(s)>-5$ by 
$$C_j(s)(\psi(\tau)+\tau^2)^{-\frac{\Re(s)+5}{2}}\tau^2.$$
In conclusion, for $\Re(s)>-2$ by Lebesgue dominated 
convergence 
$$L^<(s)(\mu)\to_{\mu\to 0} L^<(s)(0)=\frac{1}{i\pi}
\sum_j\int_{W_{\alpha_j}} A(D,s)_{<\alpha_j}.$$
But $L^<(0)(\mu)$ is constant in $\mu$, and this ends the proof of (\ref{aleg}).

\paragraph{The regularized value $\zeta(0)$.}
The zeta function of $d_t$ is treated as in Theorem
\ref{mz}. We get as before
formula (\ref{limz}) with $D$
replaced by $Q$:

\begin{eqnarray*}
t{\zeta}(d_t,s)&=&t
\Tr\left((d_t^2)^{-\frac{s}{2}}\right)\nonumber\\
&\sim_{t\to 0}& \frac{1}{2\pi}
\int_{T^*S^1} \Tr(q(d_t^2)^{-\frac{s}{2}})d\tau d\theta.
\end{eqnarray*}
The leading term in this expression is
$$\pi^{-1}\int_{T^*S^1}\Tr(Q^*Q+\tau^2)^{-s/2}
d\tau  d\theta.$$
It is clear that the 
off-diagonal terms in $q(d_t^2)$ contribute to the trace only 
modulo $t^2$. It follows that the sub-leading term
(the coefficient of $t$ in $t\overline{\zeta}(d_t,s)$) 
comes only from the two diagonal terms in $q(d_t^2)$. 
If we choose $\phi$  to be even in $\tau$ 
then these terms are odd in $\tau$, and thus 
their integral with respect to $\tau$ vanishes. This proves
(\ref{stz}) modulo (\ref{alzg}).

Let us identify $\Tr\left((Q^*Q+\tau^2)^{-\frac{s}{2}}\right)$.
As for the eta function, we split the trace according to large 
or small eigenvalues of $Q^*Q$. The large eigenvalues part 
does not depend on $\phi$ since $Q_{>\alpha_j}=D_{>\alpha_j}$; 
we compute easily 
$$\int_{T^*W_{\alpha_j}}\Tr(Q_{>\alpha_j}^*Q_{>\alpha_j}+\tau^2)^{-s/2}
d\tau d\theta=f\left(\frac{s}{2}\right)
\int_{W_{\alpha_j}}\Tr(D_{>\alpha_j}^*D_{>\alpha_j})^{-\frac{s-1}{2}}d\theta.$$
At $s=0$ we get 
$-\int_{W_{\alpha_j}}\Tr_w|D_{>\alpha_j}|d\theta$ since
$sf(s/2)_{|s=0}=-\pi$. Remark that $\Tr_w|D_{>\alpha_j}|=\Tr_w|D|$ since $\Tr_w$
vanishes on finite rank operators.
We claim next that the meromorphic function with simple poles
$$K(s):=\int_{\rz}\Tr(Q_{<\alpha_j}^*Q_{<\alpha_j}+\tau^2)^{-s/2}
d\tau$$
vanishes at $s=0$, which implies that the contribution of small eigenvalues 
is null. Indeed, let $A(\tau):=Q_{<\alpha_j}^*Q_{<\alpha_j}$. Clearly
the matrix $A$ satisfies $\|A(\tau)\|\leq C$, 
$\|A'(\tau)\|\leq C(1+\tau^2)^{-1}$ since
$\phi(\tau)$ is Schwartz. Integration by parts (for large $\Re(s)$)
shows 
\begin{eqnarray*}
K(s)&=&\int_{\rz}\tau'\Tr(A(\tau)+\tau^2)^{-s/2}
d\tau\\
&=&s/2\int_{\rz}
\Tr(\tau(2\tau+A'(\tau))(A(\tau)+\tau^2)^{-s/2-1}
d\tau\\
&=&sK(s)+s\int_\rz\Tr((\tau A'(\tau)/2 -A(\tau))
(A(\tau)+\tau^2)^{-s/2-1}d\tau.
\end{eqnarray*}
By unique continuation the identity holds for all $s\in\cz$.
For $\Re(s)\geq -1$ the last integral converges absolutely, 
so because of the $s$ factor
we deduce that $(1-s)K(s)$ is regular and vanishes at $s=0$, which is 
to say that $K(0)=0$.

Thus at $s=0$ we deduce (\ref{alzg}), (\ref{stz}) for $d_t$ 
in lieu of $\delta_t$. 
But by Lemma \ref{eeez}
we know that $\zeta(\delta_t,0)=\zeta(d_t,0)$. 
\qed

The crucial point in the above proof is reducing the analysis to a finite 
number of eigenvalues and the corresponding eigenspaces. Indeed, the use of the 
function $\phi\equiv 1$ in (\ref{q}) instead of a Schwartz
function leads to an operator in the extended adiabatic algebra if and only if 
the fibration $X\to M$ is the identity fibration $M\to M$.

Note that the sign in formula (\ref{aleg}) disappears if we 
change the definition of $\delta_t$ to be consistent with the choice of
the bounding spin structure on $S^1$.

Also, when the dimension of $N$ is even, the eta invariant is a $\pmod{2\zz}$
homotopy invariant of elliptic operators and moreover $\ind(D)=0$ because
the fibers are odd-dimensional. Thus for every $t>0$ we get:
$$e^{-i\pi\oeta(\delta_t)}=\hol(\det(D)).$$
This identity can also be deduced from the results of \cite{bf2}.

We have deliberately avoided mentioning the regularity of the eta function at 
$s=0$ since it plays no role in the proof. In light of this regularity however, 
the first part of the theorem generalizes the Witten-Bismut-Freed holonomy 
theorem to any family of elliptic first-order differential operators.

\vspace{0.5cm}\section{Closing remarks}\label{clorem}

The adiabatic algebra is a powerful tool
for studying degenerate families of operators such as $\delta_t$. Compared to previous
results in this direction, one needs to make only minimal assumptions about the family
$D$, while obtaining significant regularity in the adiabatic limit for free
from the adiabatic formalism. The smoothness of the eta and zeta functions in the
adiabatic limit is by no means obvious; it can actually fail if $D$ is non-invertible
even if $\delta_t$ is invertible for all $t$. For example, consider the 
operator on $S^1$
$$\delta_t:=\left[\begin{array}{cc}-ti{\ps}-t\alpha&0\\ 0& 
ti{\ps}+t\alpha\end{array}\right].$$ 
For general $\alpha$ and $s$, the expansion of $\eta(\delta_t,s)$ and 
$\zeta(\delta_t,s)$ at $t=0$ will
contain $\log t$ terms. The reason is that $\delta_t$ is not invertible as an adiabatic
operator.

One can consider other elliptic adiabatic operators 
constructed from the family $D$, for instance
$$P_t:=t\tna_{\ps}+D$$
in the case where $D$ is self-adjoint. 
If the fibers of $N\rightarrow S^1$ are spin, $N$ is even-dimensional 
and $D$ is the family of Dirac operators on the fibers then 
$P_t$ is the chiral Dirac operator for the metric (\ref{gt}) on $N$. We can 
apply the previous analysis 
to $\zeta(P_t,0)$. Formula (\ref{alzg}) and Theorems \ref{mz}, \ref{pzl} 
hold for $P_t$ if we divide the right-hand side by $2$.

In the particular case where $N=S^1$ is the identity fibration and 
$D=1$, Theorem \ref{mz} says that the function $t \sum_{k=-\infty}^\infty 
(t^2 k^2+1)^{-s/2}$, which is
well-defined for $(s,t)\in\{\Re(s)>1\}\times(0,\infty)$, extends as a meromorphic
function in $s\in\cz$ and smooth in $t\in[0,\infty)$. Moreover, for all $s\in\cz$,
$$\lim_{t\rightarrow 0} \ t\!\! \sum_{k=-\infty}^\infty (t^2 k^2+1)^{-s/2}= \pi^{1/2}
\frac{\Gamma\left(\frac{s-1}{2}\right)}{\Gamma\left(\frac{s}{2}\right)}.$$
This limit of meromorphic functions can be checked directly for $\Re(s)>1$, 
but seems hard to prove for general $s$ by elementary methods.
 
It is evident that
$\ind(P_t)=\zeta(P_t,0)-\zeta(P_t^*,0).$
The left-hand side is constant in $t$, in particular 
$$\ind(P_t)=\lim_{t\rightarrow 0} (\zeta(P_t,0)-\zeta(P_t^*,0)).$$
Using the adiabatic algebra methods developed in this paper, we obtain 
\begin{eqnarray*}
\zeta(P_t,0)&\sim_{t\rightarrow 0}& -\frac{1}{2t}\int_{S^1}\Tr_w(D^2)^{\frac{1}{2}}
d\theta
+\frac{1}{4}\int_{S^1}\Tr_w\left((D^2)^{-\frac{1}{2}}\tna_{\ps}
(D)\right)d\theta,\\
\zeta(P_t^*,0)&\sim_{t\rightarrow 0}& -\frac{1}{2t}\int_{S^1}\Tr_w(D^2)^{\frac{1}{2}}
d\theta
-\frac{1}{4}\int_{S^1}\Tr_w\left((D^2)^{-\frac{1}{2}}\tna_{\ps}
(D)\right)d\theta
\end{eqnarray*}
which imply
\begin{equation}\label{desf}
\ind(P_t)=\frac{1}{2}\int_{S^1}\Tr_w\left((D^2)^{-\frac{1}{2}}
\tna_{\ps}
(D)\right)d\theta.
\end{equation}
Formally, the right-hand side of (\ref{desf})
is $\frac{1}{2}\int_{S^1}\Tr_w(\tna_{\ps}(\log|D|))d\theta$,
and it is easily proved to equal
the net flow of eigenvalues through $0$ around the circle. We are therefore led to
a purely analytical proof of the identity
\begin{equation}\label{sf}
\ind(P_t)=\sfl_{S^1}(D)
\end{equation}
where $\sfl$ denotes the spectral flow \cite{aps}.

In the spirit of \cite{meleta} we can interpret the right-hand 
side of (\ref{desf}) as the integral on $S^1$ of a closed 1-form defined on $S^1$ 
(or more generally on the basis $M$ of the family $D$)  by
$$\mathrm{var}(\eta):=\Tr_w\left((D^2)^{-\frac{1}{2}}\tna(D)\right).$$
In \cite{meleta} Melrose defined the analogous quantity in the more
general setting of $1$-suspended operators, 
and proved subsequently the extension of (\ref{sf}) to that setting. Our methods
can be used to reprove Melrose's result for an elliptic adiabatic family in 
$\pa(N)$.

\bibliographystyle{plain}

\end{document}